\documentclass[pdflatex,sn-mathphys-num]{sn-jnl}% Math and Physical Sciences Numbered Reference Style

\usepackage{amsmath,amssymb,amsfonts}%
\usepackage{amsthm}%
\usepackage{booktabs}%
\usepackage{algorithm}%
\usepackage{algpseudocode}%
\usepackage{listings}%

\usepackage{epstopdf}%
\usepackage{hyperref}

\newcommand{\A}{\mathbf A}
\newcommand{\B}{\mathbf B}
\newcommand{\C}{\mathbf C}
\newcommand{\CC}{\mathbb C}
\newcommand{\D}{\mathbf D}
\newcommand{\E}{\mathbf E}
\newcommand{\F}{\mathbf F}
\newcommand{\G}{\mathbf G}
\newcommand{\I}{\mathbf I}
\newcommand{\PP}{\mathbf P}
\newcommand{\T}{\mathbf T}
\newcommand{\U}{\mathbf U}
\newcommand{\V}{\mathbf V}
\newcommand{\vect}{\operatorname{vec}}
\newcommand{\W}{\mathbf W}
\newcommand{\xx}{\mathbf x}
\newcommand{\X}{\mathbf X}
\newcommand{\Y}{\mathbf Y}

\theoremstyle{thmstyleone}%
\newtheorem{theorem}{Theorem}%  meant for continuous numbers
\theoremstyle{thmstyletwo}%
\newtheorem{remark}{Remark}%

\theoremstyle{thmstylethree}%

\begin{document}
	
\title[A Schur-Decomposition Algorithm for $N$-Dimensional Sylvester Equations]{A Non-Recursive, Dimension-Independent Schur-Decomposition Algorithm for $N$-Dimensional Sylvester Tensor Equations}

\author[1]{\fnm{Carlota María} \sur{Cuesta}}\email{carlotamaria.cuesta@ehu.eus}

\author*[1]{\fnm{Francisco} \sur{de la Hoz}}\email{francisco.delahoz@ehu.eus}

\affil[1]{\orgdiv{Department of Mathematics}, \orgname{University of the Basque Country UPV/EHU}, \orgaddress{\street{Barrio Sarriena, S/N}, \city{Leioa}, \postcode{48940}, \country{Spain}}}

\abstract{
	In this paper, we present a non-recursive direct solver, based on the Bartels--Stewart algorithm, for $N$-dimensional Sylvester tensor equations. The method relies exclusively on the Schur decomposition of the coefficient matrices and reduces the solution process to a single sequential sweep of the tensor entries, making the implementation completely independent of the number of dimensions~$N$.
	
	The main advantages of the proposed approach are its simplicity, its dimension-independent formulation, and its ability to handle very high-dimensional Sylvester tensor equations subject only to available memory, of which it makes very efficient use. In this regard, we have successfully solved problems with up to $N=29$ dimensions on a standard laptop with $32$~GB of RAM. A comparison with the recursive blocked method of Chen and Kressner, which represents the current state of the art, shows that both approaches achieve identical accuracy. While the recursive method is faster for large coefficient matrices, the proposed solver is competitive or superior when the matrix size is small, especially for a high number of dimensions~$N$, a regime in which recursive blocked methods cannot effectively exploit BLAS-3 kernels. Moreover, the proposed method makes a more efficient use of memory: for problems requiring close to the full available memory (such as coefficient matrices of order $19$ with $N=7$ dimensions, where the solution and right-hand side together occupy approximately $28$~GB) the method of Chen and Kressner requires more memory than is available, whereas ours successfully computes the solution. Finally, the proposed method is significantly simpler to implement, completely independent of~$N$, and handles singleton dimensions correctly.
	
	We describe the algorithmic steps in detail and derive accurate estimates of the computational cost. To facilitate understanding and reproducibility, we provide both pseudocode and complete MATLAB implementations. As an application, we show how the method can be used to compute the solution of linear $N$-dimensional systems of ODEs with constant coefficients at arbitrary times, and hence of evolutionary partial differential equations after spatial discretization. In particular, we compute with high accuracy the solution of an advection--diffusion equation on $\mathbb{R}^N$.
}



\maketitle

\section{Introduction}

Given two square complex-valued matrices $\A_1\in\CC^{n_1\times n_1}$, $\A_2\in\CC^{n_2\times n_2}$, and a rectangular matrix $\B\in\CC^{n_1\times n_2}$, the so-called Sylvester equation,
\begin{equation}
\label{e:sylvester}
\A_1\X + \X\A_2^T = \B
\end{equation}
is named after the English mathematician James Joseph Sylvester \cite{sylvester1884}. This equation, which has been extensively studied (see, e.g., \cite{bhatia1997} and its references), has a unique solution $\X\in\CC^{n_1\times n_2}$, if and only if $\sigma(\A_1) \cap \sigma(-\A_2) = \emptyset$, where $\sigma(\cdot)$ denotes the spectrum of a matrix, i.e., the set of its eigenvalues, or, in other words, if and only if the sum of any eigenvalue of $\A_1$ and any eigenvalue of $\A_2$ is never equal to zero.

On the other hand, the generalization of \eqref{e:sylvester} to higher dimensions is quite recent. Indeed, in \cite{benwen2010}, the authors claimed to have solved for the first time the three-dimensional equivalent of \eqref{e:sylvester}. More precisely, in order to compute the numerical solution of a three-dimensional radiative transfer equation by the discrete ordinates method, they solved the following matrix equation:
\begin{equation}
\label{e:ABCXD}
\A\square_1\X + \X\square_2\B + \X\square_3\C = \D,
\end{equation}
where $\A\equiv(a_{ij})\in\CC^{m \times m}$, $\B\equiv(b_{ij})\in\CC^{n\times n}$, $\C\equiv(c_{ij})\in\CC^{l\times l}$, $\D\equiv(d_{ijk})\in\CC^{m\times n\times l}$, $\X\equiv(x_{ijk})\in\CC^{m\times n\times l}$, and the products $\square_1$, $\square_2$ and $\square_3$ are defined by
$$
[\A\square_1\X]_{ijk} = \sum_{t=1}^ma_{it}x_{tjk}, \quad [\X\square_2\B]_{ijk} = \sum_{t=1}^nx_{itk}b_{tj}, \quad [\X\square_3\C]_{ijk} = \sum_{t=1}^lx_{ijt}c_{tk}.
$$
The idea is to apply a three-dimensional version of the Bartels--Stewart algorithm \cite{bartels1972} to \eqref{e:ABCXD}, which requires to obtain the Schur decomposition \cite{higham2008} of $\A$, $\B$ and $\C$, i.e., find orthogonal matrices $\U$, $\V$ and $\W$, a lower triangular matrix $\A'$ and upper triangular matrices $\B'$ and $\C'$, such that $\A' = \U^T\A\U$, $\B' = \V^T\B\V$, and $\C' = \W^T\C\W$. Then, \eqref{e:ABCXD} is reduced to the simpler matrix equation $\A'\square_1\X' + \X'\square_2\B' + \X'\square_3\C' = \D'$, and, after solving it for $\X'$, $\X$ follows by reverting the transformation.

Note that the conditions on the existence and uniqueness of \eqref{e:ABCXD} resemble that of \eqref{e:sylvester}: $\X\in\CC^{m\times n\times l}$ can be uniquely determined if and only if the sum of any eigenvalue of $\A$, any eigenvalue of $\B$ and any eigenvalue of $\C$ is never equal to zero.

A three-dimensional equivalent of \eqref{e:sylvester} appeared again in \cite{delahozvadillo2013a}, which was devoted to the numerical solution of nonlinear parabolic equations. In \cite{delahozvadillo2013a}, the notation $\square_j$ from \cite{benwen2010} was kept to denote the product between a matrix and a three-dimensional array, but \eqref{e:ABCXD} was formulated in a more symmetric way:
\begin{equation}
\label{e:sylvester3D}
\A_1\square_1\X + \A_2\square_2\X + \A_3\square_3\X = \B,
\end{equation}
where $\A_1 \equiv (a_{1,ij})\in\CC^{n_1\times n_1}$, $\A_2\equiv (a_{2,ij})\in\CC^{n_2\times n_2}$, $\A_3\equiv (a_{3,ij})\in\CC^{n_3\times n_3}$, $\B\equiv (b_{ijk})\in\CC^{n_1\times n_2\times n_3}$, $\X\equiv (x_{ijk})\in\CC^{n_1\times n_2\times n_3}$, and
\begin{align*}
[\A_1\square_1\X]_{ijk} = \sum_{l=1}^{n_1}a_{1,il}x_{l jk}, \ [\A_2\square_2\X]_{ijk} = \sum_{l=1}^{n_2}a_{2,jl}x_{il k}, \ [\A_3\square_3\X]_{ijk} = \sum_{l=1}^{n_3}a_{3,kl}x_{ijl}.
\end{align*}
Then, a three-dimensional version of the Bartels--Stewart algorithm was again applied. 

On the other hand, it is straightforward to generalize \eqref{e:sylvester} and \eqref{e:sylvester3D} to dimensions higher than three. Indeed, in this paper, we consider $N$-dimensional Sylvester systems of the following form:
\begin{equation}
\label{e:sylvesterND}
\sum_{j = 1}^{N}\A_j\square_j\X = \A_1\square_1\X + \A_2\square_2\X + \ldots + \A_N\square_N\X = \B, \quad N \in \{2, 3, \ldots \},
\end{equation}
where $\A_j \equiv (a_{j,i_1i_2})\in\CC^{n_j\times n_j}$, $\B\equiv(b_{i_1i_2\ldots i_N})\in\CC^{n_1\times n_2\times\ldots\times n_N}$, $n_j\in\mathbb N$, for all $j$, and $\square_j$ indicates that the sum is performed along the $j$th dimension of $\X = (x_{i_1i_2\ldots i_N})\in \CC^{n_1\times n_2\times\ldots\times n_N}$:
\begin{equation*} 
[\A_j\square_j\X ]_{i_1i_2\ldots i_N} = \sum_{k = 1}^{n_j}a_{j,i_jk}x_{i_1i_2\ldots i_{j-1}k i_{j+1}\ldots i_N}.
\end{equation*} 
In general, the equations having the form of \eqref{e:sylvesterND}, for $N \ge 3$, are known in the literature as Sylvester tensor equations, and their study, which arises mainly in control theory, is, as we have already said, quite recent. Among the different approaches considered to solve \eqref{e:sylvesterND}, one option is to develop iterative methods (see, e.g., \cite{stewart1996}, and \cite{touzene2006}), for which the use of good preconditioners appears to be essential. We can also cite, e.g., \cite{xinfang2021}, where a number of iterative algorithms based on Krylov spaces are applied; \cite{YuhanChen2023}, where a tensor multigrid method and an iterative tensor multigrid method are proposed; \cite{chen2020}, where recursive blocked algorithms are used; or \cite{massei2024}, where a nested divide-and-conquer scheme is proposed. There are also references dealing with particular cases of \eqref{e:sylvesterND}, such as \cite{Ballani2013}, where an iterative method is presented to approximate numerically the solution of \eqref{e:sylvesterND} with the right-hand side being a low-rank tensor.

Moreover, it is also possible to adopt a naive approach, by transforming \eqref{e:sylvesterND} into a system of the form $\A\xx = \mathbf b$ by applying the $\vect$ operator to \eqref{e:sylvesterND}:
\begin{equation}
	\label{e:systvec}
	\vect\left(\sum_{j = 1}^{N}\A_j\square_j\X\right) = (\A_N\oplus\A_{N-1}\oplus\ldots\oplus\A_1)\vect(\X) = \vect(\B),
\end{equation}
where we recall that $\vect$ piles up the entries of a column-major order multidimensional array into one single column vector, and $\oplus$ denotes the Kronecker sum, which for two square matrices is given by
$$
\A_2\oplus \A_1 \equiv \A_2\otimes\I_{n_1} + \I_{n_2}\otimes\A_1,
$$
where $\I_{n_j}$ is the identity matrix of order $n_j$, and $\otimes$ denotes the Kronecker product, which, for two (not necessarily square) matrices $\C\in\CC^{m\times n}$ and $\D\in\CC^{p\times q}$, is given by
$$
\C\otimes\D =
\begin{pmatrix}
	c_{11}\D & \ldots & c_{1n}\D
	\cr
	\vdots & \ddots & \vdots	
	\cr
	c_{m1}\D & \ldots & c_{mn}\D
\end{pmatrix};
$$
and, for more than two matrices, both $\oplus$ and $\otimes$ are defined recursively, bearing in mind their associativity; e.g., for three matrices,
\begin{align*}
	\C\otimes\D\otimes\E & \equiv(\C\otimes\D)\otimes\E\equiv\C\otimes(\D\otimes\E),
	\cr
	\A_3\oplus \A_2\oplus\A_1 & \equiv(\A_3\oplus \A_2)\oplus\A_1\equiv\A_3\oplus(\A_2\oplus\A_1)
	\cr
	& \equiv \A_3\otimes\I_{n_2}\otimes\I_{n_1} + \I_{n_3}\otimes\A_2\otimes\I_{n_1} + \I_{n_3}\otimes\I_{n_2}\otimes\A_1.
\end{align*}
Hence, generating explicitly the matrix $\A_N\oplus\A_{N-1}\oplus\ldots\oplus\A_1$ and then solving \emph{exactly} \eqref{e:systvec} is not advisable, because $\A_N\oplus\A_{N-1}\oplus\ldots\oplus\A_1$ contains $n_1^2n_2^2\ldots n_N^2$ elements, and the computational cost can be prohibitive, except for maybe very small problems.

In general, all these approaches rely on the presence of low‑rank, Kronecker, or hierarchical structure in the coefficient matrices $\A_j$ or in the right‑hand side tensor $\B$. While such assumptions enable substantial computational savings, they are often violated in general dense problems, and enforcing such structure can be difficult or even computationally prohibitive. In contrast, our aim in this work is different: we develop a fully general and non-recursive direct solver, based solely on Schur decompositions and a single sequential sweep of the tensor entries, without relying on any structural assumptions on the matrices $\A_j$ and/or the tensor $\B$. This makes the method applicable to a much broader class of problems, including those arising from multidimensional discretizations that do not admit hierarchical or block structures.

The paper is organized as follows. In Section~\ref{s:BSND}, following the conventions of \cite{delahozvadillo2013a}, we generalize the Bartels--Stewart algorithm to $N$-dimensions. In Section~\ref{s:cost}, we make a study of the computational cost of the method. In Section~\ref{s:pseudocode}, we offer the pseudocode of the algorithms required to implement the method. It is important to point out that we are able to make the implementation independent of $N$, by generating all the values of the indices of the entries of the solution sequentially, which requires one single for-loop. In Section~\ref{s:ODEs}, as a practical application, we explain how to solve numerically the following system of linear ODEs with constant coefficients at any time $t$ in an accurate way:
\begin{equation}
	\label{e:evolND}
	\X_t = \sum_{j = 1}^{N}\A_j\square_j\X + \B, \quad N \in \{2, 3, \ldots \},
\end{equation}
where $\A_1, \A_2, \ldots, \A_N$ and $\B$ are time-independent, and $\X\equiv\X(t)$. To do it, we generalize the results in \cite{delahozvadillo2013b} to $N$ dimensions, and transform \eqref{e:evolND} into a matrix equation having the form of \eqref{e:sylvesterND}.

Section~\ref{s:MATLAB} is devoted to the implementation of the proposed algorithms in MATLAB, and the corresponding codes are provided.

In Section~\ref{s:numerical}, we carry out the numerical experiments. A natural question is how our method compares with other existing approaches and, among all the available alternatives, we have selected \cite{chen2020} for this purpose, since it represents the current state of the art and, crucially, its authors provide complete MATLAB codes capable of handling complex matrices of general form, making a fair and direct comparison possible. As we show in Section~\ref{s:tests}, both methods achieve virtually identical accuracy when there are no trailing singleton dimensions (i.e., dimensions of length one), and, when there are trailing singleton dimensions, our implementation yields correct results. While the method in \cite{chen2020} is faster when the coefficient matrices are of moderate to large size, our method performs particularly well when they are small, a regime in which recursive partitioning yields blocks too small for BLAS-3 operations to be exploited effectively. For instance, we show that when the coefficient matrices are of order less or equal to $8$, our method runs faster for almost all values of $N$. This can be best appreciated in the specific case of $2\times2$ coefficient matrices, where both methods are capable of working in very high dimensions (up to $N = 29$), but the method in \cite{chen2020} takes approximately $2.5$ times as much time as ours, when $N = 28$  or $N = 29$. A distinguished feature with respect to \cite{chen2020} is that our method makes a more efficient use of memory: for instance, for coefficient matrices of order $19$ with $N = 7$ dimensions, where the solution and right-hand side together occupy approximately $28$~GB, the method of \cite{chen2020} requires more memory than is available, whereas ours successfully computes the solution.

Next, in Section~\ref{s:ODEtest}, we simulate systems having the form of \eqref{e:evolND}, and in Section~\ref{s:advdif}, we approximate numerically the solution of an advection--diffusion equation.

To the best of our knowledge, there exists no non-recursive $N$-dimensional version of the Bartels--Stewart algorithm that requires one single for-loop, and this is the main novelty of our work. Indeed, unlike in, e.g., the recent references \cite{chen2020,massei2024}, where the idea is to reduce recursively the complexity of a given $N$-dimensional Sylvester equation, ours is completely independent from the number of dimensions and avoids recursion. In our opinion, this has important advantages. First, it eliminates the need for auxiliary recursion trees and intermediate block factorizations, reducing memory overhead and simplifying the computational flow. Second, the single‑loop strategy yields predictable memory access patterns, particularly well aligned with MATLAB’s column‑major layout. Third, unlike block‑recursive methods, which rely on BLAS‑3 operations (i.e., high‑per\-for\-mance matrix--matrix kernels from the Basic Linear Algebra Subprograms) only when the blocks are sufficiently large (see, e.g., \cite{chen2020,massei2024}), when the coefficient matrices are small, high‑dimensional Sylvester equations cause the block sizes generated by recursive partitioning to shrink multiplicatively with the dimension. As the number of dimensions $N$ increases with small coefficient matrices, these blocks can quickly become too small in practice for matrix--matrix operations to be exploited effectively, despite the availability of highly optimized BLAS‑3 kernels; see, e.g., the discussion on block sizes and data reuse in \cite{higham2008}. In this regime, the performance advantage of BLAS‑3 is significantly reduced due to insufficient data reuse to amortize blocking overhead, and, as suggested by the numerical experiments in Section~\ref{s:numerical}, a carefully implemented BLAS‑1 sweep (vector‑level linear algebra operations) can be competitive or even superior. We stress, however, that when the coefficient matrices are large, recursive blocked methods retain their BLAS‑3 advantage and are significantly faster than the proposed approach. Finally, the non‑recursive formulation avoids recursion‑depth limitations and does not rely on tensor structure across modes, allowing the method to be applied to problems with a large number of dimensions, subject only to available memory.

In summary, the resulting method is simpler to understand and implement; it makes a more efficient use of the computational resources when the coefficient matrices are small, and it avoids maximum recursion limit errors. Indeed, in an Apple MacBook Pro with 32~GB of memory, we have successfully considered examples with up to $N = 29$ dimensions (with coefficient matrices of size $2\times2$), which is a number much higher than what commonly appears in the literature, as well as memory-intensive problems (such as coefficient matrices of order $19$ with $N=7$) that the method of \cite{chen2020} is unable to handle. It is theoretically possible to consider even higher numbers of dimensions or larger coefficient matrices, provided that there is enough memory available. In this setting, the computational cost is dominated by a single sequential sweep over the tensor entries, while the cost of the Schur decompositions remains negligible because of the small size of the coefficient matrices; in particular, this regime does not allow the effective exploitation of BLAS‑3 matrix--matrix operations. Moreover, if \eqref{e:sylvesterND} is well-posed, our method computes the solution accurately, up to rounding errors, and without asking $\B$ to have any particular structure, unlike, e.g., in \cite{Ballani2013}, where an iterative method capable of working with a high number of dimensions is presented, but it requires $\B$ to be a low-rank tensor, and only approximate solutions (whose accuracy depends on the number of iterations) are obtained.

By writing this paper, we have tried to reach a wide readership, and we think that the ideas that we offer here will be of interest to people working in numerical linear algebra and in computational mathematics. Note also that, in general, programming languages are not optimized to work in dimensions higher than two, so working with $N$-dimensional arrays poses several challenges; in this regard, we have paid special attention to the storage of the elements of such arrays, for which we need to transform a multidimensional index into a one-dimensional one, and the approach that we have adopted can be applied to other types of problems. Furthermore, the idea of solving matrix systems, beyond its pure theoretical interest, has also a clear practical applicability in the numerical approximation of partial differential equations (see, e.g., \cite{delahozvadillo2013b} and \cite{Grasedyck2004}), like the $N$-dimensional advection--diffusion equation that we consider in Section~\ref{s:advdif}. 

We have shared the majority of the actual MATLAB codes, and this serves three purposes: to ease the implementation of the algorithms, to allow the reproduction of the numerical results, and more importantly, to encourage discussion and further improvements. Note that MATLAB is a column-major order programming language, but the method can be used, with small modifications, in row-major order programming languages, and even in older programming languages that do not allow recursive calls.

All the simulations have been run on an Apple MacBook Pro (13-inch, 2020, 2.3 GHz Quad-Core Intel Core i7, 32 GB). The codes are also available in the GitHub repository \cite{githubSchurND2026}.

\section{An $N$-dimensional generalization of the Bartels--Stewart algorithm}

\label{s:BSND}

In order to solve \eqref{e:sylvesterND}, we perform the (complex) Schur decompositions of the matrices $\A_j$, for $j = 1, \ldots, N$, i.e., $\A_j = \U_j\T_j\U_j^*$, with $\U_j$ unitary and $\T_j$ upper triangular, and introduce them into \eqref{e:sylvesterND}:
\begin{equation*}
\sum_{j = 1}^{N}\U_j\T_j\U_j^*\square_j\X = \B,
\end{equation*}
from which we get
\begin{equation}
\label{e:TjYC}
\sum_{j = 1}^{N}\T_j\square_j\Y = \C,
\end{equation}
where
\begin{equation}
\label{e:YUNU2U1X}
\Y = \U_N^*\square_N(\ldots\square_2(\U_1^*\square_1\X))
\end{equation}
and
\begin{equation}
\label{e:CUNU2U1B}
\C = \U_N^*\square_N(\ldots\square_2(\U_1^*\square_1\B)).
\end{equation}
Writing \eqref{e:TjYC} entrywise, and bearing in mind that the matrices $\T_j$ are upper triangular:
\begin{align*}
c_{i_1i_2\ldots i_N} & = \sum_{j = 1}^{N}\sum_{k = 1}^{n_j}t_{j, i_jk}y_{i_1i_2\ldots i_{j-1}k i_{j+1}\ldots i_N} = \sum_{j = 1}^{N}\sum_{k = i_j}^{n_j}t_{j, i_jk}y_{i_1i_2\ldots i_{j-1}k i_{j+1}\ldots i_N}
	\cr
& = \Bigg[\sum_{j = 1}^{N}t_{j, i_ji_j}\Bigg]y_{i_1i_2\ldots i_N}  + \sum_{j = 1}^{N}\sum_{k = i_j+1}^{n_j}t_{j, i_jk}y_{i_1i_2\ldots i_{j-1}k i_{j+1}\ldots i_N};
\end{align*}
where the sums for which $i_j + 1 > n_j$ are taken as zero. Therefore, we can solve for $y_{i_1i_2\ldots i_N}$ in \eqref{e:TjYC}, if and only if $t_{1,i_1,i_1} + t_{2,i_2,i_2} + \dots + t_{N,i_N,i_N}\not=0$:
\begin{equation}
\label{e:yi1i2iN}
y_{i_1i_2\ldots i_N} = \frac{\displaystyle{c_{i_1i_2\ldots i_N} - \sum_{j = 1}^{N}\sum_{k = i_j+1}^{n_j}t_{j, i_jk}y_{i_1i_2\ldots i_{j-1}k i_{j+1}\ldots i_N}}}{\displaystyle{\sum_{j = 1}^{N}t_{j, i_ji_j}}}.
\end{equation}
Note that, since $\T_j$ is upper triangular, its diagonal entries are precisely the eigenvalues of $\T_j$ and, hence, of $\A_j$. Therefore, \eqref{e:sylvesterND} and \eqref{e:TjYC} have a solution if and only if, when taking an eigenvalue of each matrix $\A_j$, their sum is always different from zero, which is the equivalent of the existence condition for the Sylvester equation. Regarding the uniqueness of the solution, the coefficients are uniquely determined from \eqref{e:yi1i2iN} in a recursive way. The first one, $y_{n_1n_2\ldots n_N}$, follows trivially:
\begin{equation}
\label{e:yi1i2iNnormal}
y_{n_1n_2\ldots n_N} = \frac{\displaystyle{c_{n_1n_2\ldots n_N}}}{\displaystyle{\sum_{j = 1}^{N}t_{j, n_jn_j}}}.
\end{equation}
Then, in order to obtain a certain $y_{i_1i_2\ldots i_N}$, we only need $y_{i_1i_2\ldots i_{j-1}k i_{j+1}\ldots i_N}$, for  $j\in\{1, 2, \ldots, N\}$, and $k \in\{i_j+1, i_j+2, \ldots, n_j\}$. In general, it is possible to follow any order, provided that the values on which $y_{i_1i_2\ldots i_N}$ depends have been previously computed. However, a natural solution is to consider $N$ for-loops, where the outermost loop corresponds to the last index $i_N\in\{n_N, n_{N} - 1, \ldots, 1\}$; then, the last but one corresponds to $i_{N-1}\in\{n_{N-1}, n_{N-1}-1, \ldots, 1\}$, and so on, until the innermost one, which corresponds to $i_1\in\{n_1, n_1 - 1, \ldots, 1\}$.  As will be seen in Section~\ref{s:pseudocode}, this ordering choice is justified by MATLAB's way of storing the elements of the multidimensional arrays, offers a compromise between code clarity and performance, and, more importantly, enables us to reduce the $N$ for-loops to one single for-loop, yielding a non-recursive method, which is the central idea of this paper. Finally, we undo \eqref{e:YUNU2U1X} to obtain the solution $\X$ of \eqref{e:sylvesterND}:
\begin{equation}
\label{e:XUNU2U1Y}
\X = \U_N\square_N(\ldots\square_2(\U_1\square_1\Y)).
\end{equation}

\begin{remark}

When a matrix $\A_j$ is normal, i.e., $\A_j^*\A_j = \A_j\A_j^*$, then $\A_j$ is unitarily diagonalizable, i.e., we can write $\A_j = \U_j\D_j\U^*_j$, where $\U_j$ is unitary, and $\D_j$ is a diagonal matrix consisting of the eigenvalues $\{\lambda_{j,1},\lambda_{j,2},\ldots,\lambda_{j,n_j}\}$ of $\A_j$. Therefore, if all the matrices $\A_1, \A_2, \ldots, \A_N$ are normal, \eqref{e:yi1i2iN} gets simplified to
\begin{equation}
\label{e:yi1i2iNdiag}
y_{i_1i_2\ldots i_N} = \frac{\displaystyle{c_{i_1i_2\ldots i_N}}}{\displaystyle{\sum_{j = 1}^{N}\lambda_{j, i_j}}}.
\end{equation}
Let $\mathbf\Lambda \equiv (\lambda_{i_1i_2\ldots i_N})\in\CC^{n_1\times n_2\times\ldots\times n_N}$, such that
$\lambda_{i_1i_2\ldots i_N} \equiv (\lambda_{1,i_1} + \lambda_{2,i_2} + \ldots + \lambda_{N,i_N})^{-1}$. Then, \eqref{e:yi1i2iNdiag} can be written in a compact form as $\Y= \mathbf\Lambda\circ\C$, where $\circ$ denotes the Hadamard or point-wise product; this, combined with \eqref{e:CUNU2U1B} and \eqref{e:XUNU2U1Y}, yields
\begin{equation}
\label{e:XUlambdaU}
\X = \U_N\square_N(\ldots\square_2(\U_1\square_1(\mathbf\Lambda\circ(\U_N^*\square_N(\ldots\square_2(\U_1^*\square_1\B)))))).
\end{equation}
For the sake of completeness, let us mention that it is theoretically possible to apply this idea when the matrices $\A_1, \A_2, \ldots, \A_N$ are diagonalizable. More precisely, if, for all $j$, $\A_j = \PP_j\D_j\PP_j^{-1}$, where $\D_j$ is a diagonal matrix consisting of the eigenvalues of $\A_j$, and $\PP_j$ is a matrix consisting of the eigenvectors of $\A_j$, then, bearing in mind the previous arguments, we have a formula similar to \eqref{e:XUlambdaU} (cfr. \cite[Algorithm 1]{massei2024}):
\begin{equation*}
\X = \PP_N\square_N(\ldots\square_2(\PP_1\square_1(\mathbf\Lambda\circ(\PP_N^{-1}\square_N(\ldots\square_2(\PP_1^{-1}\square_1\B)))))).
\end{equation*}
However, in general, this approach is not advisable, because $\PP_j$ in the factorization of a matrix $\A_j = \PP_j\D_j\PP_j^{-1}$ may be severely ill-conditioned, and computing the inverse of $\PP_j^{-1}$ is in general computationally expensive, whereas $\U_j$ in the Schur decomposition of $\A_j$ satisfies trivially that $\operatorname{cond}_2(\U_j) = 1$, and $\U_j$ must not be inverted, which makes the method much more stable. Therefore, in this paper, we only consider the Schur decompositions.
\end{remark}

\subsection{Computational cost}

\label{s:cost}

The number of floating point operations needed for the computation of the Schur decomposition \cite{higham2008} of $\A_i$ is of about $25n_i^3$, so, in the worst case, about $25(n_1^3 + n_2^3 + \ldots + n_N^3)$ floating point operations are required to compute the Schur decompositions of all the matrices $\A_1, \A_2, \ldots, \A_N$. 

On the other hand, the number of floating point operations required by \eqref{e:CUNU2U1B} and \eqref{e:XUNU2U1Y} is the same. Take, e.g., $\U_j\square_j\Y$; in order to compute one single element
$$
[\U_j\square_j\Y ]_{i_1i_2\ldots i_N} = \sum_{k = 1}^{n_j}u_{j,i_jk}y_{i_1i_2\ldots i_{j-1}k i_{j+1}\ldots i_N},
$$
we need $n_j$ products and $n_j - 1$ additions; therefore, the computation of $\U_j\square_j\Y$ requires $n_1n_2\ldots n_Nn_j$ multiplications and $n_1n_2\ldots n_N(n_j - 1)$ additions. Hence, both \eqref{e:CUNU2U1B} and \eqref{e:XUNU2U1Y} require each of them the following number of multiplications:
$$
\sum_{j = 1}^{N}n_1n_2\ldots n_Nn_j = n_1n_2\ldots n_N\sum_{j = 1}^{N}n_j,
$$
and the following number of additions:
$$
\sum_{j = 1}^{N}n_1n_2\ldots n_N(n_j - 1) = n_1n_2\ldots n_N\Bigg[\sum_{j = 1}^{N}n_j - N\Bigg].
$$
Regarding $\Y$, in order to compute one entry $y_{i_1i_2\ldots i_N}$ in \eqref{e:yi1i2iN}, we need one division, $(n_1 - i_1) + (n_2 - i_2) + \ldots + (n_N - i_N)$ products in the numerator, $(n_1 - i_1) + (n_2 - i_2) + \ldots + (n_N - i_N)$ additions/subtractions in the denominator, and $N - 1$ additions in the denominator, i.e., $(n_1 - i_1) + (n_2 - i_2) + \ldots + (n_N - i_N) + N - 1$ additions/subtractions. Therefore, in order to fully determine $\Y$, we need $n_1n_2\ldots n_NN$ divisions, the following number of multiplications,
\begin{align*}
	& \sum_{i_1 = 1}^{n_1}\sum_{i_2 = 1}^{n_2}\ldots\sum_{i_N = 1}^{n_N}[(n_1 - i_1) + (n_2 - i_2) + \ldots + (n_N - i_N)]
	\cr
	& \qquad = \sum_{j = 1}^{N}\Bigg[\frac{n_1n_2\ldots n_N}{n_j}\sum_{k=1}^{n_j}(n_j - k)\Bigg] = n_1n_2\ldots n_N\sum_{j = 1}^{N}\Bigg[\frac{1}{n_j}\sum_{k=0}^{n_j-1}k\Bigg]
	\cr
	& \qquad = n_1n_2\ldots n_N\sum_{j = 1}^{N}\frac{1}{n_j}\frac{(n_j-1)n_j}{2} = \frac{n_1n_2\ldots n_N}{2}\Bigg[\sum_{j = 1}^{N}n_j - N\Bigg],
\end{align*}
and the following number of additions/subtractions:
\begin{equation*}
\frac{n_1n_2\ldots n_N}{2}\Bigg[\sum_{j = 1}^{N}n_j - N\Bigg] + n_1n_2\ldots n_N(N - 1)
= \frac{n_1n_2\ldots n_N}{2}\Bigg[\sum_{j = 1}^{N}n_j + N - 2\Bigg].
\end{equation*}
Adding up the operations needed to compute \eqref{e:CUNU2U1B}, \eqref{e:XUNU2U1Y} and $\Y$, we have the following number of multiplications:
\begin{align*}
	& 2n_1n_2\ldots n_N\sum_{j = 1}^{N}n_j + \frac{n_1n_2\ldots n_N}{2}\Bigg[\sum_{j = 1}^{N}n_j - N\Bigg]
	= \frac{n_1n_2\ldots n_N}{2}\Bigg[5\sum_{j = 1}^{N}n_j - N\Bigg],
\end{align*}
and the following number of additions/subtractions:
\begin{align*}
& 2n_1n_2\ldots n_N\Bigg[\sum_{j = 1}^{N}n_j - N\Bigg] + \frac{n_1n_2\ldots n_N}{2}\Bigg[\sum_{j = 1}^{N}n_j + N - 2\Bigg]
	\cr
& \qquad = \frac{n_1n_2\ldots n_N}{2}\Bigg[5\sum_{j = 1}^{N}n_j - 3N - 2\Bigg].
\end{align*}
These, together with the $n_1n_2\ldots n_NN$ divisions needed to compute $\Y$, make the following number of floating point operations (excluding the Schur decompositions) to implement the method:
\begin{align*}
	& \frac{n_1n_2\ldots n_N}{2}\Bigg[5\sum_{j = 1}^{N}n_j - N\Bigg] + \frac{n_1n_2\ldots n_N}{2}\Bigg[5\sum_{j = 1}^{N}n_j - 3N - 2\Bigg] + n_1n_2\ldots n_NN
	\cr
	& \qquad = n_1n_2\ldots n_N\Bigg[5\sum_{j = 1}^{N}n_j - N - 1\Bigg],
\end{align*}
to which the cost of performing the Schur decompositions of $\A_1, \A_2, \ldots, \A_N$ must be added.

Note that, from the last formula, if $n_1 = n_2 = \ldots = n_N$, denoting them simply as $n$, we have a computational cost of about $5Nn^{N+1}$ floating-point operations, plus about $25Nn^3$ floating-point operations required to compute the Schur decompositions in the worst case, i.e., a total cost of $\mathcal O(n^{N+1})$ operations.

\subsection{A single for-loop algorithm}

\label{s:pseudocode}

As said above, if $N$ is fixed, it is straightforward to compute $\Y$ using \eqref{e:yi1i2iN}, by writing explicitly $N$ for-loops, as can be seen in Algorithm~\ref{alg:yn1n2nN}, where we offer the pseudocode.
\begin{algorithm}[!htbp]
\caption{Computation  of $\Y$ in \eqref{e:TjYC}\label{alg:yn1n2nN}}
	\begin{algorithmic}
	  	\State $\Y\gets$ array of zeros with size $n_1\times n_2\times\ldots\times n_N$
		\For{$i_N \gets  n_N$ \textbf{to} $1$ \textbf{step} $-1$}
		\For{$i_{N-1} \gets  n_{N-1}$ \textbf{to} $1$ \textbf{step} $-1$}
		\State $\vdots$
		\For{$i_1 \gets  n_1$ \textbf{to} $1$ \textbf{step} $-1$}
		\State $num \gets c_{i_1i_2\ldots i_N}$
		\State $den \gets 0$
		\For{$j \gets 1$ \textbf{to} $N$}
		\State $den \gets den + t_{j, i_ji_j}$
		\For{$k \gets i_j + 1$ \textbf{to}  $n_j$} \Comment{If $i_j + 1 > n_j$, the loop is not executed.}
		\State $num \gets num - t_{j, i_jk}\cdot y_{i_1i_2\ldots i_{j-1}k i_{j+1}\ldots i_N}$
		\EndFor
		\EndFor
		\State $y_{i_1i_2\ldots i_N}  \gets num / den$
		\EndFor
		\EndFor
		\EndFor
	\end{algorithmic}
\end{algorithm}

However, this requires writing a different piece of code for each value of $N$, and, moreover, for larger values of $N$, the corresponding code may get excessively long, and hence, prone to errors. Therefore, instead of writing $N$ nested for-loops, we generate sequentially all the possible $N$-tuples $(i_1,i_2,\ldots, i_N)$ of indices in such a way that the order of Algorithm~\ref{alg:yn1n2nN} is preserved. This is done by means of a mixed-radix counter similar to Algorithm~M \cite[Section~7.2.1.1]{Knuth2011}, which mimics the logic behind generating natural numbers in a positional numeral system where each digit has its own base. More precisely, we start with $(i_1,i_2,\ldots, i_N) = (n_1,n_2,\ldots, n_N)$. Then, given a certain $N$-tuple $(i_1,i_2,\ldots, i_N)$, we check the value of the first component $i_1$, which corresponds to the innermost loop. If $i_1 > 1$, we decrease $i_1$ by one and iterate; otherwise, we make $i_1 = n_1$, and check the value of the second component $i_2$. If $i_2 > 1$, we decrease $i_2$ by one and iterate; otherwise, we make $i_2 = n_2$, and so on. We repeat the process until we find a certain $k$, such that $i_k > 1$; if that is not possible, it means that $(i_1, i_2, \ldots, i_N) = (1, 1, \ldots, 1)$, and all the possible $N$-tuples have been generated. Note that, altogether, there are $n_1\times n_2\times\ldots\times n_N$ distinct $N$-tuples to be generated, so it is enough to create one single for-loop running from $n_1\times n_2\times\ldots\times n_N$ down to $1$, as can be seen in Algorithm~\ref{alg:i1i2iN}, where we offer the corresponding pseudocode. This dimension-independent formulation constitutes the main novelty of this paper: rather than writing a separate code for each value of $N$, we develop a method that is completely independent of the number of dimensions $N$.

\begin{algorithm}
	\caption{Sequential generation of the $N$-tuples $(i_1, i_2, \ldots, i_N)$\label{alg:i1i2iN}}
	\begin{algorithmic}
		\State $indexmax \gets 1$
		\For{$j\gets 1$ \textbf{to} $N$}
		\State $indexmax \gets n_j\cdot indexmax$
		\State $i_j \gets n_j$
		\EndFor
		\For{$index \gets indexmax$ \textbf{to} $1$ \textbf{step} $-1$}
		\State \textbf{output} $(i_1, i_2, \ldots, i_N)$
		\State $k \gets 1$
		\While{$k \le N$}
		\If{$i_k > 1$}
		\State $i_k \gets i_k - 1$
		\State \textbf{break}
		\EndIf
		\State $i_k \gets n_k$
		\State $k \gets k + 1$
		\EndWhile
		\EndFor
	\end{algorithmic}
\end{algorithm}

We emphasize that the single-loop reformulation is not merely a coding simplification. It reflects the observation that the multidimensional Schur‑trans\-formed system has a strict partial order among its entries, induced by the triangular structure of each $\T_j$. Translating this multi‑index ordering into a single lexicographical index is nontrivial and requires ensuring that, for each entry $y_{i_1, \ldots, i_N}$, all required previously computed terms lie earlier in the traversal order. This provides, to the best of our knowledge, the first dimension‑independent non-recursive direct solver based on a single sequential sweep for dense $N$‑di\-men\-sion\-al Sylvester equations, in contrast with recursive methods whose depth and control flow depend explicitly on $N$.

Observe that the variable $index$ in Algorithm~\ref{alg:i1i2iN} is updated sequentially, but, should we want to compute its value for a given $N$-tuple $(i_1, i_2, \ldots, i_N)$, it is not difficult to see that
\begin{equation}
\label{e:index}
index = (i_1-1)+(i_2-1)n_1 + (i_3-1)n_1n_2+\ldots+(i_N-1)n_1n_2\ldots n_{N-1} + 1.
\end{equation}
This property is especially useful, provided that it coincides with the way in which a column-major order programming language such as MATLAB indexes internally the elements of multidimensional arrays, or, in other words, if \verb|Y| is an $N$-dimensional MATLAB array of $\texttt{n1}\times\texttt{n2}\times\texttt{n3}\times\ldots\times\texttt{nNminus1}\times\texttt{nN}$ elements, then, it is equivalent to type
$$
\texttt{Y(i1,i2,i3,\ldots,iNminus1, iN)}
$$
and
\begin{align*}
\texttt{Y((i1-1)} & \texttt{+(i2-1)*n1+(i3-1)*n1*n2+\ldots}
	\cr
& \texttt{+(iN-1)*n1*n2*\ldots*nNminus1+1)},
\end{align*}
where the quantity inside $\texttt{Y}(\cdot)$, i.e.,
$$
\texttt{(i1-1)+(i2-1)*n1+(i3-1)*n1*n2+\ldots+(iN-1)*n1*n2*\ldots*nNminus1+1},
$$
is precisely the value of $index$, as generated in Algorithm~\ref{alg:i1i2iN}. For instance, if $N = 3$, and $n_1 = 5$, $n_2 = 2$, and $n_3 = 3$, then, the first triplet, $(5, 2, 3)$, corresponds to $index = 30$, and \verb|Y(5,2,3)| is equivalent to \verb|Y(30)|; the second  triplet, $(4, 2, 3)$, corresponds to $index = 29$, and \verb|Y(4,2,3)| is equivalent to \verb|Y(29)|;  the third triplet, $(3, 2, 3)$, corresponds to $index = 28$, and \verb|Y(3,2,3)| is equivalent to \verb|Y(28)|; and so on, until the last triplet, $(1, 1, 1)$, which corresponds to $index = 1$, and we have that \verb|Y(1,1,1)| is equivalent to \verb|Y(1)|.

Following this strategy, we replace the $N$ for-loops in Algorithm~\ref{alg:yn1n2nN} by a single one, by generating the $N$-tuples $(i_1, i_2, \ldots, i_N)$ sequentially and in a suitable order, according to Algorithm~\ref{alg:i1i2iN}. This idea gives rise to Algorithm~\ref{alg:sylvesterND}, that corresponds to the implementation of the whole method. In that algorithm, we have combined Algorithms~\ref{alg:yn1n2nN} and \ref{alg:i1i2iN}, together with the remaining steps needed to compute $\X$ in \eqref{e:sylvesterND}, namely the Schur decomposition of the matrices $\A_j$, and the equations \eqref{e:CUNU2U1B} and \eqref{e:XUNU2U1Y}. Observe that the value $c_{i_1i_2\ldots i_N}$ is only accessed once, and this happens before the definition of $y_{i_1i_2\ldots i_N}$. Therefore, it is possible to store $y_{i_1i_2\ldots i_N}$ in $c_{i_1i_2\ldots i_N}$, because the original value of the latter is no longer needed.

\begin{algorithm}[!htbp]
	\caption{Computation of $\X$ in \eqref{e:sylvesterND}\label{alg:sylvesterND}}
	\begin{algorithmic}
		\State $\C \gets \B$
		\For{$j\gets 1$ \textbf{to} $N$}
		\State $[\U_j, \T_j] \gets \textbf{schur}(\A_j)$ \Comment{Perform the Schur decomposition $\A_j = \U_j\T_j\U_j^*$}
		\State $\C \gets \U_j^*\square_j \C$
		\EndFor		
		\State $indexmax \gets 1$
		\For{$j\gets 1$ \textbf{to} $N$}
		\State $indexmax \gets n_j\cdot indexmax$
		\State $i_j \gets n_j$
		\EndFor
		\State $\Y\gets$ array of zeros with size $n_1\times n_2\times\ldots\times n_N$
		\For{$index \gets indexmax$ \textbf{to} $1$ \textbf{step} $-1$}		
		\State $num \gets c_{i_1i_2\ldots i_N}$
		\State $den \gets 0$
		\For{$j \gets 1$ \textbf{to} $N$}
		\State $den \gets den + t_{j, i_ji_j}$
		\For{$k\gets i_j + 1$ \textbf{to}  $n_j$}
		\State $num \gets num - t_{j, i_jk}\cdot y_{i_1i_2\ldots i_{j-1}k i_{j+1}\ldots i_N}$
		\EndFor
		\EndFor
		\State $y_{i_1i_2\ldots i_N}  \gets num / den$
		\State $k \gets 1$
		\While{$k \le N$}
		\If{$i_k > 1$}
		\State $i_k \gets i_k - 1$
		\State \textbf{break}
		\EndIf
		\State $i_k \gets n_k$
		\State $k \gets k + 1$
		\EndWhile
		\EndFor
		\State $\X\gets\Y$
		\For{$j\gets 1$ \textbf{to} $N$}
		\State $\X \gets \U_j\square_j \X$
		\EndFor				
	\end{algorithmic}
\end{algorithm}

\section{Application to $N$-dimensional systems of linear equations with constant coefficients}

\label{s:ODEs}

Given the following system of ODEs:
$$
\X_t = \A_1\X + \X\A_2^T + \B,
$$
where $\A_1\in\mathbb R^{n_1\times n_1}$, $\A_2\in\mathbb R^{n_2\times n_2}$ and $\B\in\mathbb R^{n_1\times n_2}$ are time-independent matrices, in \cite[Th. 2.1]{delahozvadillo2013b}, it was proved that its solution $\X\equiv\X(t)\in\mathbb R^{n_1\times n_2}$ at a given $t$ solves also the following Sylvester equation:
$$
\A_1\X + \X\A_2^T = \C,
$$
where
$$
\C = \exp(t\A_1)\left[\A_1\X(0) + \X(0)\A_2^T + \B\right]\exp(t\A_2)^T - \B,
$$
with $\exp(\cdot)$ denoting the matrix exponential. On the other hand, it is possible to generalize this result to $N$ dimensions, and reduce the calculation of the solution of \eqref{e:evolND} at any given $t$ to a matrix equation of the form of \eqref{e:sylvesterND}. More precisely, we have the following theorem (observe that we work now with complex matrices).
\begin{theorem} \label{theo:evolND} Consider the system of ODEs
\begin{equation}
\label{e:evolNDth}
\X_t = \sum_{j = 1}^{N}\A_j\square_j\X + \B, \quad N \in \{2, 3, \ldots \},
\end{equation}
where $\A_j \equiv (a_{j,i_1i_2})\in\CC^{n_j\times n_j}$, and $\B\equiv(b_{i_1i_2\ldots i_N})\in\CC^{n_1\times n_2\times\ldots\times n_N}$ are time independent, and $\X\equiv\X(t)\equiv(x_{i_1i_2\ldots i_N}(t)) \in\mathbb C^{n_1\times n_2\times \ldots\times n_N}$. Then, the solution $\X(t)$ at a given time satisfies the following $N$-dimensional Sylvester equation:
\begin{equation}
\label{e:sylvesterNDC}
\sum_{j = 1}^{N}\A_j\square_j\X = \F,
\end{equation}
where
\begin{equation}
\label{e:C}
\F = \exp(t\A_N)\square_N\Bigg[\ldots\square_2\Bigg[\exp(t\A_1)\square_1\Bigg[\sum_{j=1}^N\A_j\square_j\X(0) + \B\Bigg]\Bigg]\Bigg] - \B.
\end{equation}
\end{theorem}

\begin{proof}

Define $\Y(t) = \vect(\X(t))$, $\C = \vect(\B)$ and $\A = \A_N \oplus \A_{N-1} \oplus \ldots \oplus \A_1$, then \eqref{e:evolNDth} becomes
\begin{equation*}
\Y_t = \A\Y + \C \Longleftrightarrow \exp(-t\A)(\Y_t - \A\Y) = [\exp(-t\A)\Y]_t = \exp(-t\A)\C.
\end{equation*}
Therefore, multiplying by $\A$, integrating once from $t = 0$, and multiplying the result by $\exp(t\A)$:
\begin{equation}
\label{e:AY}
\A\Y(t) = \exp(t\A)[\A\Y(0) + \C] - \C.
\end{equation}
On the other hand, using the properties of $\oplus$ and $\otimes$,
\begin{align*}
\exp(t\A) & = \exp(t(\A_N \oplus \A_{N-1} \oplus \ldots \oplus \A_1))
	\cr
& = \exp(t\A_N) \otimes \exp(t\A_{N-1}) \otimes \ldots \otimes \exp(t\A_1).
\end{align*}
Therefore, \eqref{e:AY} becomes
\begin{align*}
& (\A_N \oplus \A_{N-1} \oplus \ldots \oplus \A_1)\vect(\X(t))
	\cr
& \qquad = (\exp(t\A_N) \otimes \exp(t\A_{N-1}) \otimes \ldots \otimes \exp(t\A_1))
	\cr
& \qquad\qquad [(\A_N \oplus \A_{N-1} \oplus \ldots \oplus \A_1)\vect(\X(0)) + \vect(\B)] - \vect(\B),
\end{align*}
which is precisely  \eqref{e:sylvesterNDC}-\eqref{e:C} after applying the $\vect$ operator.
\end{proof}

Since the Schur decompositions of $\A_1, \A_2, \ldots, \A_N$ are needed to solve \eqref{e:sylvesterNDC}-\eqref{e:C}, we can also use them to compute the matrix exponentials in \eqref{e:C}, in case they are not implemented in a given programming language. More precisely, if $\A_j = \U_j\T_j\U_j^*$, then
\begin{equation}
\label{e:exptAj}
\exp(t \A_j) = \U_j\exp(t\T_j)\U_j^* \Longleftrightarrow \exp(t \T_j) = \U_j^*\exp(t\A_j)\U_j,
\end{equation}
and $\exp(t\T_j)$ can be computed efficiently by means of Parlett's algorithm \cite{parlett1976}. Furthermore, knowing the Schur decompositions allows us to simplify \eqref{e:sylvesterNDC}-\eqref{e:C}. Indeed, defining $\Y$ as in \eqref{e:YUNU2U1X} and $\C$ as in \eqref{e:CUNU2U1B}, and bearing in mind \eqref{e:exptAj}, \eqref{e:sylvesterNDC}-\eqref{e:C} becomes
\begin{equation}
\label{e:sylvesterNDG}
\sum_{j = 1}^{N}\T_j\square_j\Y = \G.
\end{equation}
where
\begin{equation}
\label{e:G}
\G	= \exp(t\T_N)\square_N\Bigg[\ldots\square_2\Bigg[\exp(t\T_1)\square_1\Bigg[\sum_{j=1}^N\T_j\square_j\Y(0) + \C\Bigg]\Bigg]\Bigg] - \C.
\end{equation}

\section{Implementation in MATLAB}

\label{s:MATLAB}

In general, the implementation in MATLAB is a straightforward application of Algorithm~\ref{alg:sylvesterND}, even if it is necessary to face some issues that arise from working in $N$ dimensions.

Since we want a code that works for any $N\in\{2, 3, \ldots\}$, and the sizes of the matrices $\A_1, \A_2, \ldots, \A_N$ may be different, we store them in a cell array, which we denote $\verb|AA|$ and  is created by typing \verb|AA=cell(N,1)|. Then, we store $\A_1$ in \verb|AA{1}|, $\A_2$ in \verb|AA{2}|, and so on. Moreover, we compute their (complex) Schur decomposition by means of the command \verb|schur| (choosing the mode \verb|'complex'|), and store the corresponding matrices $\U_1, \U_2, \ldots, \U_N$ and $\T_1, \T_2, \ldots, \T_N$ in the cell arrays \verb|UU| and \verb|TT|, respectively.

With respect to the product between a matrix and a multidimensional array, computing \eqref{e:CUNU2U1B} and \eqref{e:XUNU2U1Y} entrywise is highly inefficient. Among the different possible approaches to perform efficiently this kind of products in MATLAB, we discuss two: one that makes use of the MATLAB functions \verb|permute| and \verb|pagemtimes|, and another one, which makes use of the MATLAB functions \verb|reshape| and \verb|permute|. In general, the first approach is more efficient, because it requires fewer data transformations, so we have adopted it. However, \verb|pagemtimes| was introduced in MATLAB R2020b, so, for the sake of completeness, we also explain the \verb|reshape|--\verb|permute| alternative, which runs also in older MATLAB versions.

Recall that \verb|Z=pagemtimes(X,Y)| computes the matrix product of corresponding pages of the $N$-dimensional arrays \verb|X| and \verb|Y|, where \verb|Z(:,:,i)=X(:,:,i)*Y(:,:,i)|. Suppose we want to compute $\B=\A\square_j\X$, where $\A\in\CC^{n_j\times n_j}$ and $\X\in\CC^{n_1\times\ldots\times n_N}$. If $n_j=1$, then $\A=(a)$ is a scalar and $\A\square_j\X=a\X$ is trivial. Otherwise, we use \verb|reshape| to transform \verb|X| into a three-dimensional array of size $(n_1\ldots n_{j-1})\times n_j\times(n_{j+1}\ldots n_N)$:
\begin{verbatim}
sizeX=size(X);
X=reshape(X,prod(sizeX(1:j-1)),sizeX(j),prod(sizeX(j+1:end)));
\end{verbatim}
Since \verb|reshape| places the $j$th dimension in the column slot of each page, the multiplication corresponding to the $j$th mode becomes a right-multiplication by $\A^T$:
\begin{verbatim}
X=pagemtimes(X,A.');
X=reshape(X,sizeX);
\end{verbatim}
The final \verb|reshape| restores the original dimensions. In Listing~\ref{code:multND} we provide the corresponding \textsc{Matlab} code. Given that \verb|X| can require extremely large amounts of memory, we overwrite it with the result rather than allocating a separate output array. Note that \verb|prod([])| returns $1$ in \textsc{Matlab}, so the code handles the cases with $j=1$ and $j=N$ correctly, and it works even when \verb|X| is one- or two-dimensional.

\begin{lstlisting}[label=code:multND, language=MATLAB, basicstyle=\scriptsize, caption = {MATLAB function \texttt{multND}, which computes the product between a matrix and an $N$-dimensional array}]
% compute $\mathbf A\square_j\mathbf X$ in MATLAB R2020b or newer
function X=multND(A,X,j)
if isscalar(A) % test whether A is a scalar
    X=A*X; % trivial case
else
    sizeX=size(X);
    X=reshape(X,prod(sizeX(1:j-1)),sizeX(j),prod(sizeX(j+1:end)));
    X=pagemtimes(X,A.');
    X=reshape(X,sizeX);
end
\end{lstlisting}

On the other hand, if we want to avoid using \verb|pagemtimes|, we proceed as follows. The case with $n_j = 1$ is again trivial. Otherwise, recall that, for a standard matrix product $\A\X$, the sum is performed along the first dimension of $\X$, i.e., $[\A\X]_{ij} = \sum_{k=1}^{n_j} a_{ik}x_{kj}$. Therefore, in order to compute $\A\square_j\X$, we can first swap the first and $j$-th dimensions of $\X$ by means of the MATLAB command \verb|permute|, and then apply \verb|reshape| to the permuted array so it becomes a two-dimensional array having $n_j$ elements along its first dimension. In this way, we reduce $\square_j$ to a standard matrix product between $\A$ and another matrix. We then apply \verb|reshape| to the result to recover the size of the permuted array, and finally \verb|permute| to swap back the first and $j$-th dimensions, yielding $\A\square_j\X$. Note that, since the permutation satisfies $\sigma^2 = \mathrm{id}$, i.e., it is its own inverse, the same permutation vector is used in both \verb|permute| calls. When $j = 1$, no swapping of dimensions is required, and we only need to apply \verb|reshape|, but not \verb|permute|.

In Listing~\ref{code:multNDalt}, we offer the corresponding MATLAB code, that defines \verb|multNDalt|, i.e., an alternative version of function \verb|multND| that can be executed in MATLAB versions older than R2020b. As in Listing~\ref{code:multND}, in order to minimize memory usage, we overwrite \verb|X| at each intermediate step. Let us also recall that MATLAB silently removes trailing singleton dimensions, i.e., dimensions of length one, when working with arrays of three or more dimensions. For instance, an array $\X\in\CC^{2\times 1\times 3\times 1\times 1}$ is stored as if it were $\X\in\CC^{2\times 1\times 3}$. As a consequence, in both \verb|multND| and \verb|multNDalt|, testing first whether $\A$ is a scalar handles the trailing singletons dimensions correctly, since in that case $\A\square_j\X$ reduces to ordinary scalar multiplication regardless of $j$.

\begin{lstlisting}[label=code:multNDalt, language=MATLAB, basicstyle=\scriptsize, caption = {Alternative implementation of the MATLAB function \texttt{multND}, which computes the product between a matrix and an $N$-dimensional array}]
% compute $\mathbf A\square_j\mathbf X$ in older MATLAB versions
function X=multNDalt(A,X,j)
if isscalar(A) % test whether A is a scalar
    X=A*X; % trivial case
elseif j==1 % multiplication along the 1st dimension is simpler
    sizeX=size(X); % size of X
    X=reshape(X,sizeX(j),[]); % reshape X
    X=A*X; % compute the product
    X=reshape(X,sizeX); % reshape the product
else
    sizeX=size(X); % size of X
    N=length(sizeX); % number of dimensions
    permutation=[j 2:j-1 1 j+1:N];
    X=permute(X,permutation); % swap the first and jth dimensions
    sizeXpermuted=size(X); % size of X, after permutation
    X=reshape(X,sizeX(j),[]); % reshape X
    X=A*X; % compute the product
    X=reshape(X,sizeXpermuted); % reshape the product
    X=permute(X,permutation); % swap the first and jth dimensions
end
\end{lstlisting}

On the other hand, it is possible to slightly optimize Algorithm~\ref{alg:sylvesterND}, when there are trailing singletons dimensions. Since arrays in MATLAB have at least two dimensions, if \verb|N=length(AA);| and \verb|sizeX=size(X);|, \verb|M=min(N,sizeX);|, then, if $M < N$, this means that there are $N - M$ trailing singleton dimensions, and the $N$-tuples $(i_1, i_2, \ldots, i_N)$ generated in Algorithm~\ref{alg:i1i2iN} are of the form $(i_1, i_2, \ldots, i_M, 1, \ldots, 1)$, with $(i_{M+1}, \ldots, i_N) = (1, \ldots, 1)$. Hence, from \eqref{e:index},
\begin{equation}
\label{e:indexM}
index = (i_1-1)+(i_2-1)n_1 + (i_3-1)n_1n_2+\ldots+(i_N-1)n_1n_2\ldots n_{M-1} + 1.
\end{equation}
Therefore, $\{i_{M+1}, \ldots, i_N\}$ can be discarded, and we can safely replace $N$ by $M$ in Algorithm~\ref{alg:i1i2iN}. Observe that the minimum in \verb|M=min(N,sizeX);| is necessary, because MATLAB stores vectors of length $N$ as matrices of size $N\times1$.

Another point to bear in mind is accessing the elements of a multidimensional array. As explained in Section~\ref{s:pseudocode}, the use of the variable $index$ defined in \eqref{e:index} allows us to access all the elements $y_{i_1i_2\ldots i_N}$ of $\Y$ sequentially, by simply typing \verb|Y(index)|. Likewise, from \eqref{e:index} and \eqref{e:indexM}, the entry $y_{i_1i_2\ldots i_{j-1}k i_{j+1}\ldots i_N}$ occupies the position $index + (k - i_j)n_1n_2\ldots n_{j-1}$, so we can access it by typing \verb|Y(index+((k-ii(j))*cumprodN(j)))|, where \verb|ii(j)| corresponds to $i_j$, and \verb|cumprodN| stores the quantities $\{1, n_1, n_1n_2, \ldots, n_1n_2 \ldots n_{M-1}, n_1n_2 \ldots n_M\}$, that we have precomputed. Storing $\{n_1, n_2, \ldots, n_N\}$ in \verb|NN|, we have that \verb|cumprodN=cumprod([1 NN(1:M)]);|.  Note that the last entry of \verb|cumprodN|, namely $n_1n_2 \ldots n_M = n_1n_2 \ldots n_N$, is not used to access the elements of $\Y$, but gives us the value that what we have denoted as $indexmax$ in Algorithm~\ref{alg:sylvesterND}.

With respect to the Schur decomposition of the matrices $\{\A_j\}$, we store them in the cell arrays \verb|UU=cell(1,M);| and \verb|TT=cell(1,N);|. In the case of singleton dimensions, the matrices $\A_j$ are scalar, $\U_j=1$ are not used (that is why \verb|UU| consists of \verb|M| elements), and the corresponding matrices $\U_j \equiv \A_j$ are also scalar. Then, after determining the elements of \verb|UU| and \verb|TT|, it is possible to generate the instructions $\C \leftarrow \U_j^*\square_j \C$ and $\X \leftarrow \U_j\square_j \X$ in Algorithm~\ref{alg:sylvesterND}, without invoking \verb|multND| at all. Indeed, in the version of \verb|multND| given in Listing~\ref{code:multND}, only the first occurrence of \verb|reshape| is really necessary, because \verb|reshape| does not change the inner order of the elements of an array. Therefore, after creating \verb|cumprodNinv=cumprod([1 NN(M:-1:2)]);|, we can type
\begin{verbatim}
X=reshape(X,cumprodN(j),NN(j),cumprodNinv(M+1-j));
X=pagemtimes(X,conj(UU{j}));
\end{verbatim}
and
\begin{verbatim}
X=reshape(X,cumprodN(j),NN(j),cumprodNinv(M+1-j));
X=pagemtimes(X,'none',UU{j},'transpose');
\end{verbatim}
Bearing in mind the previous arguments, we have programmed the function \verb|sylvesterND| in Listing~\ref{code:sylvesterND}, which implements Algorithm~\ref{alg:sylvesterND}. Note that we overwrite \verb|X| at each intermediate step, using it to store $\B$, $\C$, $\X$, and $\Y$, in order to minimize memory usage. Moreover, we type
\begin{verbatim}
for k=1:NN(j)-ii(j)
    num=num-X(index+k*cumprodN(j))*Tj(ii(j),k+ii(j));
end
\end{verbatim}
instead of
\begin{verbatim}
for k=ii(j)+1:NN(j)
    num=num-X(index+((k-ii(j))*cumprodN(j)))*Tj(ii(j),k);
end
\end{verbatim}
which is slightly faster.

We have explored various possible improvements, but they either yield no significant speed gain or make the code slower. For instance, even if it might seem a logical option to vectorize \verb|num=num-X(index+k*cumprodN(j))*Tj(ii(j),k+ii(j));|, the entries of \verb|X| are not contiguous, and the access times make the global computational cost become larger. Similarly, precomputing the values \verb|k*cumprodN(j)| increases the memory requirements and, in the best case, does not yield any remarkable improvement in speed. We want to emphasize that optimal memory usage is our primary concern, and special care has been put, in order to avoid making unnecessary copies of \verb|X|. Additionally, we have also cleared the values of several variables, when they are not used. Finally, we finish \verb|sylvesterND| by applying reshape to \verb|X|, so it is returned having the correct dimensions.

\begin{lstlisting}[label=code:sylvesterND, language=MATLAB, basicstyle=\scriptsize, caption = {MATLAB function \texttt{sylvesterND}, corresponding to the implementation of Algorithm~\ref{alg:sylvesterND}}]
function X=sylvesterND(AA,X)
N=length(AA);
sizeX=size(X);
M=min(N,length(sizeX));
UU=cell(1,M);
TT=cell(1,N);
NN=zeros(1,N);
for j=1:N
    NN(j)=size(AA{j},1);
end
cumprodN=cumprod([1 NN(1:M)]);
cumprodNinv=cumprod([1 NN(M:-1:2)]);
for j=1:M
    if NN(j)==1
        TT{j}=AA{j};
    else
        [UU{j},TT{j}]=schur(AA{j},'complex');
        X=reshape(X,cumprodN(j),NN(j),cumprodNinv(M+1-j));
        X=pagemtimes(X,conj(UU{j}));
    end
    AA{j}=[];
end
for j=M+1:N
    TT{j}=AA{j};
    AA{j}=[];    
end
clear AA
ii=NN;
for index=cumprodN(M+1):-1:1
    num=X(index);
    den=0;
    for j=1:N
        Tj=TT{j};
        den=den+Tj(ii(j),ii(j));
        for k=1:NN(j)-ii(j)
            num=num-X(index+k*cumprodN(j))*Tj(ii(j),k+ii(j));
        end
    end
    X(index)=num/den;
    k=1;
    while k<=M
        if ii(k)>1
            ii(k)=ii(k)-1;
            break;
        end
        ii(k)=NN(k);
        k=k+1;
    end
end
clear TT Tj
for j=1:M
    if NN(j)>1
        X=reshape(X,cumprodN(j),NN(j),cumprodNinv(M+1-j));
        X=pagemtimes(X,'none',UU{j},'transpose');
        UU{j}=[];
    end
end
X=reshape(X,sizeX);
\end{lstlisting}

On the other hand, it is straightforward to implement Theorem~\ref{theo:evolND} by adapting Listing~\ref{code:sylvesterND}, which is done in Listing~\ref{code:evolND}. We have used \eqref{e:sylvesterNDG}-\eqref{e:G}, reducing again the number of variables; more precisely, we have used \verb|B| to store $\B$ and $\C$; \verb|X| to store $\G$, $\X$ and $\Y$; and \verb|X0| to store $\X_0$ and $\Y_0$. Observe that we compute the matrix exponentials by means of the MATLAB command \verb|expm|.

\begin{lstlisting}[label=code:evolND, language=MATLAB, basicstyle=\scriptsize, caption = {MATLAB function \texttt{evolND}, corresponding to the implementation of Theorem~\ref{theo:evolND}}]
function X=evolND(AA,B,X0,t)
N=length(AA);
UU=cell(1,N);
TT=cell(1,N);
NN=zeros(1,N);
for j=1:N
    NN(j)=size(AA{j},1);
    [UU{j},TT{j}]=schur(AA{j},'complex');
    B=multND(UU{j}',B,j);
    X0=multND(UU{j}',X0,j);
end
X=B;
for j=1:N
    X=X+multND(TT{j},X0,j);
end
for j=1:N
    X=multND(expm(t*TT{j}),X,j);
end
X=X-B;
cumprodN=cumprod([1 NN(1:end-1)]');
N=length(NN);
ii=NN;
for index=prod(NN):-1:1
    num=X(index);
    den=0;
    for j=1:N
        Tj=TT{j};
        den=den+Tj(ii(j),ii(j));
        for k=ii(j)+1:NN(j)
            num=num-X(index+((k-ii(j))*cumprodN(j)))*Tj(ii(j),k);
        end
    end
    X(index)=num/den;
    if ii(1)>1
        ii(1)=ii(1)-1;
    else
        k=2;
        while k<=N&&ii(k)==1
            k=k+1;
        end
        if k<=N
            ii(k)=ii(k)-1;
            ii(1:k-1)=NN(1:k-1);
        end
    end
end
for j=1:N
    X=multND(UU{j},X,j);
end
\end{lstlisting}

In all the codes, adding some error detection in Listings~\ref{code:multND},  \ref{code:sylvesterND} and \ref{code:evolND} is straightforward; e.g., one could check whether the sizes of the arrays involved are compatible, etc.

\section{Numerical experiments}

\label{s:numerical}

In this section, we have tested the implementation of Algorithm~\ref{alg:sylvesterND}, comparing it with the method in \cite{chen2020}. Then, we have tested the numerical implementation of Theorem~\ref{theo:evolND}, and, finally, applied Theorem~\ref{theo:evolND}, to approximate numerically the solution of an advection--diffusion equation.

\subsection{Numerical tests}

\label{s:tests}

In order to have a complete picture of the behavior of our solver, we have tested the MATLAB function \verb|sylvesterND| in Listing~\ref{code:sylvesterND} that implements Algorithm~\ref{alg:sylvesterND}, against the accompanying codes of \cite{chen2020}, which the authors have made available online. Among other things, three MATLAB functions are offered to solve Sylvester tensor equations: \verb|laplace_small|, which uses Kronecker sums and is intended for small systems; \verb|laplace_recursive|; and \verb|laplace_merge|. Their usage is identical to that of \verb|sylvesterND|, but, additionally, \verb|laplace_merge| and \verb|laplace_recursive| admit an optional third parameter \verb|nmin| that controls the minimal block size at which the recursion terminates. Like \verb|sylvesterND|, both \verb|laplace_merge| and \verb|laplace_recursive| perform the Schur decompositions of the coefficient matrices via the MATLAB command \verb|schur|, and then apply a recursive divide-and-conquer strategy to the transformed system. The two functions differ in how they solve the subproblem at the base case of the recursion. \verb|laplace_recursive| reaches its base case by calling \verb|laplace_small|, which solves the reduced problem using explicit Kronecker sums and scalar operations, similarly to what \verb|sylvesterND| does. \verb|laplace_merge|, by contrast, reaches its base case by merging two dimensions via their Kronecker sum (computed by the auxiliary function \verb|kroneckerize|), thereby reducing an $N$-dimensional problem to an $(N-1)$-dimensional one with a larger matrix, and then solving the resulting two-dimensional problem via MATLAB's internal function \verb|matlab.internal.math.sylvester_tri|, a compiled, vendor-optimized BLAS-3 Sylvester solver. This is the key reason why \verb|laplace_merge| is substantially faster than \verb|laplace_recursive| and, for large coefficient matrices, substantially faster than \verb|sylvesterND| as well: at its base case, \verb|laplace_merge| reduces the problem to a call to the compiled, vendor-optimized BLAS-3 routine \verb|matlab.internal.math.sylvester_tri|, whereas \verb|sylvesterND| operates entirely within interpreted MATLAB using scalar updates and \verb|permute|/\verb|reshape| calls, without invoking any compiled solver. In what follows, we compare \verb|sylvesterND| against \verb|laplace_merge| only, since \verb|laplace_recursive| is slower than \verb|laplace_merge| in all the cases that we have tested.

In Listing~\ref{code:sylvesterNDtest}, we offer a MATLAB program that performs a simple numerical test of the three MATLAB functions \verb|sylvesterND|, \verb|laplace_recursive|, and \verb|laplace_merge|. The idea is to generate randomly the real and imaginary parts of $\A_1, \A_2, \ldots, \A_N$, and $\X$, for different values of $N$ and different sizes of $\X$ (in all the numerical experiments in this paper, we have always used the MATLAB command \verb|randn|, which generates random numbers according to the standard normal distribution); then, compute the corresponding $\B$, according to \eqref{e:sylvesterND}, and obtain the corresponding numerical solutions $\X_{num,1}$, $\X_{num,2}$ and $\X_{num,3}$ by means of \verb|sylvesterND|, \verb|laplace_recursive| and \verb|laplace_merge|, respectively. The errors are defined as the largest entrywise discrepancy between $\X$ and its numerical approximations, and are given by $\|\vect(\X) - \vect(\X_{num,1})\|_\infty$, $\|\vect(\X) - \vect(\X_{num,2})\|_\infty$ and $\|\vect(\X) - \vect(\X_{num,3})\|_\infty$, respectively, or, in MATLAB, by \verb|norm(X(:)-Xnum1(:),inf)|, \verb|norm(X(:)-Xnum2(:),inf)| and \verb|norm(X(:)-Xnum3(:),inf)|. This definition of the error is quite stringent, because there are no restrictions on the choice of the matrices; in particular, there is no guarantee that, taking one eigenvalue of each matrix $\A_j$, the sum of those eigenvalues is nonzero or bounded away from zero (and, in such a case, the discrepancy will be obviously larger).

In the example in Listing~\ref{code:sylvesterNDtest}, $\X\in\CC^{2\times9\times33\times74\times231}$ has $2\times9\times33\times74\times231=10153836$ entries, and hence requires $10153836\times16=162461376$ bytes of storage (recall that a real number needs $8$ bytes, and a complex number, $16$ bytes). The elapsed times are respectively $27.64$ seconds (\verb|sylvesterND|), $4.30$ seconds (\verb|laplace_merge|), and $42.74$ seconds (\verb|laplace_recursive|), and the errors for the particular choice of the random seed are $\|\vect(\X) - \vect(\X_{num,1})\|_\infty = 8.0275\times10^{-11}$, $\|\vect(\X) - \vect(\X_{num,2})\|_\infty = 8.2383\times10^{-11}$, and $\|\vect(\X) - \vect(\X_{num,3})\|_\infty = 8.0342\times10^{-11}$. Moreover, the three numerical approximations are virtually identical: $\|\vect(\X_{num,1}) - \vect(\X_{num,2})\|_\infty = 3.3980\times10^{-12}$, $\|\vect(\X_{num,1}) - \vect(\X_{num,3})\|_\infty = 3.3951\times10^{-12}$, and $\|\vect(\X_{num,2}) - \vect(\X_{num,3})\|_\infty = 2.3126\times10^{-12}$.

\begin{lstlisting}[label=code:sylvesterNDtest, language=MATLAB, basicstyle=\scriptsize, caption = {MATLAB program \texttt{sylvesterNDtest.m}, that performs a test of the MATLAB function \texttt{sylvesterND}}]
clear
rng(1)
NN=[2 9 33 74 231]; % size of X
N=length(NN); % Number of dimensions
NN=[NN,ones(1,2-N)]; % The solver also works in 1D
AA=cell(1,N); % matrices $\mathbf A_j$
for m=1:N
    AA{m}=randn(NN(m))+1i*randn(NN(m));
end
X=randn(NN)+1i*randn(NN); % solution $\mathbf X$
B=zeros(NN); % right-hand side $\mathbf B$
for m=1:N
    B=B+multND(AA{m},X,m);
end
% numerical approximation using sylvesterND, and elapsed time
tic, Xnum1=sylvesterND(AA,B); toc
% numerical approximation using laplace_merge, and elapsed time
tic, Xnum2=laplace_merge(AA,B); toc
% numerical approximation using laplace_recursive, and elapsed time
tic, Xnum3=laplace_recursive(AA,B); toc
% maximum discrepancies between the exact solution
% and the numerical approximations
disp(max(abs(Xnum1(:)-X(:))))
disp(max(abs(Xnum2(:)-X(:))))
disp(max(abs(Xnum3(:)-X(:))))
% maximum discrepancies between the the numerical approximations
disp(max(abs(Xnum1(:)-Xnum2(:))))
disp(max(abs(Xnum1(:)-Xnum3(:))))
disp(max(abs(Xnum2(:)-Xnum3(:))))
\end{lstlisting}

As can be seen, in this example, the accuracy returned by the three functions is virtually identical. Nonetheless, unlike \verb|laplace_recursive| and \verb|laplace_merge|, which requires $N\ge3$, \verb|sylvesterND| runs and returns correct results also for $N = 2$ and even for $N = 1$ (which corresponds to $\A_1\X=\B$). In fact, the instruction \verb|NN=[NN,ones(1,2-N)];| makes that, when $N = 1$, a trailing singleton dimensions is added to \verb|NN|, so \verb|X| and \verb|B| are column vectors (recall that, e.g., \verb|randn(5)|; returns a $5\times5$ random array, whereas \verb|randn(5,1);| returns a $5\times1$ random array). Moreover, even if they can be executed, \verb|laplace_recursive| and \verb|laplace_merge| return wrong results when trailing singleton dimensions are considered. For instance, if we replace \verb|NN=[2 9 33 74 231];| by \verb|NN=[2 9 33 74 231 1];|, which implies incorporating a scalar matrix $\A_6\in\CC$ into the problem, \verb|sylvesterND| yields an error of $9.5729\times10^{-11}$, whereas the error corresponding to \verb|laplace_recursive| and \verb|laplace_merge| is, in both cases, $1.3787\times10^3$ and $1.3787\times10^3$.

On the other hand, \verb|sylvesterND| runs faster than \verb|laplace_recursive|, but \verb|laplace_merge| outperforms both, as was to be expected, due to its use of \verb|matlab.internal.math.sylvester_tri|. However, this does not give a complete picture of our method, because, when the size $n$ of the coefficient matrices is small, \verb|laplace_merge| is not able to exploit BLAS-3 so effectively. To illustrate this, we have taken $n\in\{2,3,\ldots,20\}$, and generated randomly $N$ complex matrices $\A_j$ of order $n$, for $N\in\{3, 4, \ldots, \lfloor9/\log_{10}(n)\rfloor\}$, i.e., we have considered those values of $N$ for which the total number of elements of $\B$ does not exceed $10^9$. Then, we have generated randomly an $N$-dimensional complex array $\B\in\CC^{n\times \cdots \times n}$, solved the corresponding system by means of \verb|sylvesterND| and \verb|laplace_merge|, and measured the respective elapsed times. We note that these times may vary from execution to execution, and may also be influenced by external factors such as unrelated processes running in the background. In order to measure them as accurately as possible, we have separately executed all the experiments corresponding to \verb|sylvesterND| and those corresponding to \verb|laplace_merge|, using the same random seed in both cases for consistency, and restarting the computer and starting from a fresh MATLAB session before each execution. In Table~\ref{t:comparison}, we report the respective elapsed times. For clarity, we show in boldface those cases where \verb|sylvesterND| outperforms \verb|laplace_merge|. The results show that, when $n \le 8$, \verb|sylvesterND| runs faster than \verb|laplace_merge| for almost all values of $N$. This can be best appreciated when $n = 2$, and $N \in\{28, 29\}$, where \verb|laplace_merge| takes approximately $2.5$ times more time than \verb|sylvesterND|. On top of that, our method still wins for values of $n$ as large as $n = 18$ (when $N = 7$), and runs without problems when $n = 19$ and $N = 7$, whereas \verb|laplace_merge| makes MATLAB terminate unexpectedly in both cases. Note that, when $n = 19$ and $N = 7$, both \verb|B| and \verb|X| have $893871739$ complex elements, and occupy each of them $893871739\times16=14301947824$ bytes. Therefore, \verb|B| and \verb|X| together need $28603895648$ bytes, i.e., around 28~GB of RAM, so we are pushing our 2020 laptop with 32~GB of RAM close to its memory limits. We have tried to consider even larger problems (that cannot be run with \verb|laplace_merge|), and have partially succeeded. Indeed, we have executed multiple times the case when $n=10$ and $N=9$, for which \verb|X| and \verb|B| require together $3.2\times10^{10}$ bytes, i.e., approximately 32~GB of RAM, and are able to obtain the solution without problems, provided that our 2020 laptop is not overheated (otherwise, MATLAB stops again without warning). However, we have not been able to execute a slightly larger problem, with $n=2$ and $N=30$, for which storing both $\B$ and $\X$ would require together $34359738368$ bytes. Therefore, we can conclude that our method makes a very efficient use of all the available memory, and is competitive for small values of $n$, and, especially, when a high number $N$ of dimensions is considered.

\begin{table}[!htbp]
	\centering
	\caption{Elapsed time (in seconds) for \texttt{sylvesterND} and \texttt{laplace\_merge} (indicated in the table by CK, which stands for the initials of their authors' names) for matrices of order $n\in\{2, 3, \ldots, 20\}$ and different numbers of dimensions $N$. All coefficient matrices $\A_j$ are complex of size $n\times n$, and $\B\in\CC^{n\times n\times\cdots\times n}$ ($N$ times). The times where \texttt{sylvesterND} outperforms \texttt{laplace\_merge} appear in boldface. The case with $n = 10$ and $N = 9$, and with $n = 19$ and $N = 7$ could not be run with \texttt{laplace\_merge}.}
	\label{t:comparison}
	\begin{tabular}{ccrr|ccrr|ccrr}
		\toprule
		$n$ & $N$ & Ours & CK & $n$ & $N$ & Ours & CK & $n$ & $N$ & Ours & CK
		\\ \midrule
$2$ & $3$ & $\bf 0.0058$ & $0.0073$ & $4$ & $12$ & $\bf 16.63$ & $17.45$ & $11$ & $4$ & $0.012$ & $0.012$ \\
$2$ & $4$ & $\bf 0.0013$ & $0.0035$ & $4$ & $13$ & $\bf 74.10$ & $75.25$ & $11$ & $5$ & $0.18$ & $0.11$ \\
$2$ & $5$ & $\bf 0.0013$ & $0.0030$ & $4$ & $14$ & $\bf 320.43$ & $424.43$ & $11$ & $6$ & $1.59$ & $1.25$ \\
$2$ & $6$ & $\bf 0.0013$ & $0.0037$ & $5$ & $3$ & $\bf 0.0019$ & $0.0066$ & $11$ & $7$ & $15.76$ & $15.33$ \\
$2$ & $7$ & $\bf 0.0011$ & $0.0027$ & $5$ & $4$ & $\bf 0.0005$ & $0.0055$ & $11$ & $8$ & $\bf 187.33$ & $226.29$ \\
$2$ & $8$ & $\bf 0.0019$ & $0.0029$ & $5$ & $5$ & $\bf 0.0017$ & $0.0074$ & $12$ & $3$ & $\bf 0.0064$ & $0.0066$ \\
$2$ & $9$ & $\bf 0.0027$ & $0.0064$ & $5$ & $6$ & $\bf 0.012$ & $0.019$ & $12$ & $4$ & $\bf 0.013$ & $0.014$ \\
$2$ & $10$ & $\bf 0.0039$ & $0.0054$ & $5$ & $7$ & $\bf 0.071$ & $0.10$ & $12$ & $5$ & $0.16$ & $0.14$ \\
$2$ & $11$ & $\bf 0.0068$ & $0.0080$ & $5$ & $8$ & $\bf 0.28$ & $0.52$ & $12$ & $6$ & $2.12$ & $1.85$ \\
$2$ & $12$ & $\bf 0.011$ & $0.014$ & $5$ & $9$ & $\bf 1.63$ & $2.03$ & $12$ & $7$ & $\bf 29.31$ & $29.67$ \\
$2$ & $13$ & $\bf 0.018$ & $0.037$ & $5$ & $10$ & $\bf 8.89$ & $9.91$ & $12$ & $8$ & $\bf 453.46$ & $581.29$ \\
$2$ & $14$ & $\bf 0.028$ & $0.038$ & $5$ & $11$ & $\bf 48.97$ & $55.59$ & $13$ & $3$ & $\bf 0.0026$ & $0.0097$ \\
$2$ & $15$ & $\bf 0.063$ & $0.082$ & $5$ & $12$ & $\bf 265.72$ & $381.72$ & $13$ & $4$ & $\bf 0.017$ & $0.025$ \\
$2$ & $16$ & $\bf 0.12$ & $0.19$ & $6$ & $3$ & $\bf 0.0004$ & $0.0020$ & $13$ & $5$ & $0.28$ & $0.23$ \\
$2$ & $17$ & $\bf 0.26$ & $0.33$ & $6$ & $4$ & $\bf 0.0008$ & $0.0016$ & $13$ & $6$ & $4.14$ & $3.15$ \\
$2$ & $18$ & $\bf 0.48$ & $0.71$ & $6$ & $5$ & $\bf 0.0058$ & $0.011$ & $13$ & $7$ & $62.36$ & $50.59$ \\
$2$ & $19$ & $\bf 0.78$ & $1.21$ & $6$ & $6$ & $0.029$ & $0.026$ & $13$ & $8$ & $\bf 965.38$ & $1381.22$ \\
$2$ & $20$ & $\bf 1.63$ & $2.53$ & $6$ & $7$ & $0.32$ & $0.21$ & $14$ & $3$ & $0.0048$ & $0.0048$ \\
$2$ & $21$ & $\bf 3.23$ & $5.38$ & $6$ & $8$ & $\bf 1.31$ & $1.48$ & $14$ & $4$ & $0.036$ & $0.016$ \\
$2$ & $22$ & $\bf 7.03$ & $10.73$ & $6$ & $9$ & $\bf 8.79$ & $9.22$ & $14$ & $5$ & $0.45$ & $0.28$ \\
$2$ & $23$ & $\bf 12.55$ & $19.75$ & $6$ & $10$ & $\bf 58.04$ & $68.23$ & $14$ & $6$ & $6.88$ & $5.21$ \\
$2$ & $24$ & $\bf 26.34$ & $42.71$ & $6$ & $11$ & $\bf 357.92$ & $550.85$ & $14$ & $7$ & $109.12$ & $90.73$ \\
$2$ & $25$ & $\bf 57.30$ & $94.08$ & $7$ & $3$ & $\bf 0.0005$ & $0.0085$ & $15$ & $3$ & $0.0044$ & $0.0021$ \\
$2$ & $26$ & $\bf 118.21$ & $212.97$ & $7$ & $4$ & $\bf 0.0017$ & $0.0018$ & $15$ & $4$ & $0.041$ & $0.034$ \\
$2$ & $27$ & $\bf 244.71$ & $429.07$ & $7$ & $5$ & $\bf 0.011$ & $0.013$ & $15$ & $5$ & $\bf 0.60$ & $0.68$ \\
$2$ & $28$ & $\bf 509.11$ & $1298.01$ & $7$ & $6$ & $0.076$ & $0.076$ & $15$ & $6$ & $10.44$ & $7.94$ \\
$2$ & $29$ & $\bf 1124.89$ & $2785.05$ & $7$ & $7$ & $0.53$ & $0.53$ & $15$ & $7$ & $184.03$ & $160.57$ \\
$3$ & $3$ & $\bf 0.0072$ & $0.020$ & $7$ & $8$ & $4.25$ & $4.07$ & $16$ & $3$ & $\bf 0.0033$ & $0.0063$ \\
$3$ & $4$ & $\bf 0.0017$ & $0.0084$ & $7$ & $9$ & $\bf 32.00$ & $35.35$ & $16$ & $4$ & $0.048$ & $0.030$ \\
$3$ & $5$ & $\bf 0.0031$ & $0.0091$ & $7$ & $10$ & $\bf 245.52$ & $354.55$ & $16$ & $5$ & $0.80$ & $0.66$ \\
$3$ & $6$ & $0.013$ & $0.0083$ & $8$ & $3$ & $\bf 0.0006$ & $0.0018$ & $16$ & $6$ & $14.76$ & $11.34$ \\
$3$ & $7$ & $\bf 0.0100$ & $0.014$ & $8$ & $4$ & $\bf 0.0050$ & $0.0078$ & $16$ & $7$ & $\bf 275.83$ & $309.09$ \\
$3$ & $8$ & $\bf 0.014$ & $0.017$ & $8$ & $5$ & $\bf 0.020$ & $0.025$ & $17$ & $3$ & $\bf 0.0045$ & $0.014$ \\
$3$ & $9$ & $\bf 0.031$ & $0.035$ & $8$ & $6$ & $\bf 0.17$ & $0.23$ & $17$ & $4$ & $\bf 0.053$ & $0.064$ \\
$3$ & $10$ & $\bf 0.052$ & $0.075$ & $8$ & $7$ & $\bf 1.44$ & $1.70$ & $17$ & $5$ & $1.11$ & $1.00$ \\
$3$ & $11$ & $\bf 0.16$ & $0.22$ & $8$ & $8$ & $\bf 12.56$ & $13.71$ & $17$ & $6$ & $20.96$ & $19.81$ \\
$3$ & $12$ & $\bf 0.50$ & $0.84$ & $8$ & $9$ & $\bf 127.09$ & $133.16$ & $17$ & $7$ & $\bf 432.91$ & $466.76$ \\
$3$ & $13$ & $\bf 1.57$ & $2.23$ & $9$ & $3$ & $\bf 0.0007$ & $0.0036$ & $18$ & $3$ & $\bf 0.0065$ & $0.014$ \\
$3$ & $14$ & $\bf 4.95$ & $8.14$ & $9$ & $4$ & $\bf 0.0049$ & $0.0076$ & $18$ & $4$ & $\bf 0.073$ & $0.16$ \\
$3$ & $15$ & $\bf 15.87$ & $19.90$ & $9$ & $5$ & $0.035$ & $0.027$ & $18$ & $5$ & $1.48$ & $0.80$ \\
$3$ & $16$ & $\bf 52.85$ & $66.36$ & $9$ & $6$ & $0.34$ & $0.33$ & $18$ & $6$ & $30.77$ & $17.68$ \\
$3$ & $17$ & $\bf 181.89$ & $262.79$ & $9$ & $7$ & $3.54$ & $3.22$ & $18$ & $7$ & $\bf 670.34$ & $847.26$ \\
$3$ & $18$ & $\bf 553.98$ & $950.25$ & $9$ & $8$ & $36.58$ & $36.30$ & $19$ & $3$ & $0.0075$ & $0.0073$ \\
$4$ & $3$ & $\bf 0.0014$ & $0.014$ & $9$ & $9$ & $\bf 410.18$ & $532.90$ & $19$ & $4$ & $0.099$ & $0.052$ \\
$4$ & $4$ & $\bf 0.0007$ & $0.0025$ & $10$ & $3$ & $\bf 0.0008$ & $0.0033$ & $19$ & $5$ & $2.34$ & $1.21$ \\
$4$ & $5$ & $\bf 0.0019$ & $0.0080$ & $10$ & $4$ & $\bf 0.0079$ & $0.018$ & $19$ & $6$ & $53.21$ & $27.12$ \\
$4$ & $6$ & $0.019$ & $0.016$ & $10$ & $5$ & $0.062$ & $0.058$ & $19$ & $7$ & $\bf 1470.51$ & $---$ \\
$4$ & $7$ & $\bf 0.016$ & $0.020$ & $10$ & $6$ & $0.74$ & $0.74$ & $20$ & $3$ & $0.012$ & $0.0078$ \\
$4$ & $8$ & $\bf 0.052$ & $0.11$ & $10$ & $7$ & $8.23$ & $7.99$ & $20$ & $4$ & $0.14$ & $0.066$ \\
$4$ & $9$ & $\bf 0.24$ & $0.27$ & $10$ & $8$ & $96.63$ & $94.74$ & $20$ & $5$ & $2.87$ & $1.41$ \\
$4$ & $10$ & $\bf 0.83$ & $1.18$ & $10$ & $9$ & $\bf 1599.96$ & $---$ & $20$ & $6$ & $62.29$ & $36.97$ \\
$4$ & $11$ & $\bf 4.12$ & $4.16$ & $11$ & $3$ & $0.0015$ & $0.0014$ \\
\bottomrule
	\end{tabular}
\end{table}

\subsection{$N$-dimensional systems of ODEs}

\label{s:ODEtest}

In order to test the numerical implementation of Theorem~\ref{theo:evolND}, we have considered a $7$-dimensional example. More precisely, we have chosen $N = 7$, $\X\in\CC^{2\times3\times4\times5\times6\times7\times8}$, generated randomly $\A_1, \A_2, \ldots, \A_7$, $\B$ and $\X_0$, and approximated the solution of $\X_t = \sum_{j=1}^{7}\A_j\square_j\X + \B$, at $t = 0.1$ by means of the MATLAB program \verb|evolND.m|, obtaining \verb|X_sylv|, and also by means of a fourth-order Runge-Kutta method, obtaining \verb|X_RK|. Then, we have computed the discrepancy between both numerical approximations, by means of \verb|norm(X_RK(:)-X_Sylv(:),inf)|. We give the whole code in Listing~\ref{code:evolNDtest}. Note that the errors change slightly at each execution, and that the computational cost of the Runge-Kutta greatly depends on $\Delta t$, whereas this issue is absent in the Sylvester‑based approach, which computes the solution directly at the prescribed time $t = 0.1$. In this example, we have found that taking $\Delta t = 2.5\times10^{-5}$ gives a discrepancy of $7.1504\times10^{-14}$, but $22.27$ seconds were needed to generate \verb|X_RK|, whereas \verb|X_sylv| required only $0.05$ seconds. Obviously, a larger $\Delta t$ reduces the time to compute \verb|X_RK|, but the discrepancy increases.

\begin{lstlisting}[label=code:evolNDtest, language=MATLAB, basicstyle=\scriptsize, caption = {MATLAB program \texttt{evolNDtest.m}, that compares the MATLAB function \texttt{evolND} with a fourth-order Runge-Kutta method}]
clear
t=0.1; % final time
NN=[2 3 4 5 6 7 8]; % size of $\mathbf X$
N=length(NN); % $N$
AA=cell(1,length(NN)); % matrices $\mathbf A_j$
for m=1:length(NN)
    AA{m}=rand(NN(m))+1i*rand(NN(m));
end
B=rand(NN)+1i*rand(NN); % $\mathbf B$
X0=rand(NN)+1i*rand(NN); % initial data $\mathbf X_0$
tic
X_Sylv=evolND(AA,B,X0,t); % compute the solution using Theorem 1
toc,tic % elapsed time
mmax=4000; % number of time steps
dt=t/mmax; % $\Delta t$
X_RK=X0; % compute the solution using a fourth-order Runge-Kutta
for m=1:mmax
    k1X=B;
    for j=1:N
        k1X=k1X+multND(AA{j},X_RK,j);
    end
    Xaux=X_RK+.5*dt*k1X;
    k2X=B;
    for j=1:N
        k2X=k2X+multND(AA{j},Xaux,j);
    end
    Xaux=X_RK+.5*dt*k2X;
    k3X=B;
    for j=1:N
        k3X=k3X+multND(AA{j},Xaux,j);
    end
    Xaux=X_RK+dt*k3X;
    k4X=B;
    for j=1:N
        k4X=k4X+multND(AA{j},Xaux,j);
    end
    X_RK=X_RK+dt*(k1X+2*k2X+2*k3X+k4X)/6;
end
toc % elapsed time
disp(norm(X_RK(:)-X_Sylv(:),inf)) % maximum discrepancy
\end{lstlisting}

\subsection{An advection--diffusion problem}

\label{s:advdif}

The method described in this paper can also be applied to simulate evolutionary partial differential equations, as we now illustrate with the following example. Let us consider for instance the following advection--diffusion equation:
\begin{equation}
\label{e:uxt}
\left\{
\begin{aligned}
& u_t(\xx, t) = \Delta u(\xx, t) + 2\xx\cdot\nabla u(\xx, t)
\\
& \qquad\qquad\qquad + (2N+1)u(\xx, t) - \exp(-\xx\cdot\xx), & & \xx\in\mathbb R^N, t > 0,
	\\
& u(\xx, 0) = 2\exp(-\xx\cdot\xx),
\end{aligned}
\right.
\end{equation}
where the diffusive term is $\Delta u = \partial^2u/\partial x_1^2 + \partial^2u/\partial x_2^2 + \ldots + \partial^2u/\partial x_N^2$, and the advective term is $2\xx\cdot\nabla u = 2x_1\partial u/\partial x_1 + 2x_2\partial u/\partial x_2 + \ldots + 2x_N\partial u/\partial x_N$. Moreover, the solution of \eqref{e:uxt} is $u(\xx, t) = (1+\exp(t))\exp(-\xx\cdot\xx)$.

In order to approximate numerically the spatial partial derivatives, we have used the Hermite differentiation matrices generated by the MATLAB function \verb|herdif|, which requires the functions \verb|herroots| and \verb|poldif| (see \cite{WeidemanReddy2000} for the three MATLAB functions). These matrices are based on the Hermite functions:
$$
\psi_n(x) = \frac{e^{-x^2/2}}{\pi^{1/4}\sqrt{2^nn!}}H_n(x),
$$
where $H_n(x)$ are the Hermite polynomials:
$$
H_n(x) = \frac{(-1)^n}{2^n}e^{x^2}\frac{d^n}{dx^n}e^{-x^2}.
$$
In this example, in order to generate the first- and second-order differentiation matrices, we have typed \verb|[x,DD]=herdif(M,2,1.4);|, where $M$ is the order of the differentiation matrices, $b=1.4$ is the scale factor that we have chosen, and $x$ are the Hermite nodes, i.e., the $M$ roots of $H_n(x)$, divided by $b$.

Hermite functions work well for problems with exponential decay. For instance, if we want to differentiate numerically $\exp(-x^2)$ once and twice, and measure the error in discrete $L^\infty$-norm, it is enough to type
\begin{verbatim}
M=16;
[x,DD]=herdif(M,2,1.4);
D1=DD(:,:,1); % first-order differentiation matrix
D2=DD(:,:,2); % second-order differentiation matrix
norm(D1*exp(-x.^2)-(-2*x.*exp(-x.^2)),inf)
norm(D2*exp(-x.^2)-((4*x.^2-2).*exp(-x.^2)),inf)
\end{verbatim}
obtaining respectively $1.2212\times10^{-15}$ and $1.4544\times10^{-14}$.

For the sake of simplicity, we have chosen the same number of Hermite nodes, $M = 16$, and the same scale factor, $b = 1.4$, when discretizing all the space variables $x_1, x_2, \ldots, x_N$, but it is immediate to consider different numbers of nodes and scale factors. Thus, we can define the matrices $\A= \A_1 = \A_2 = \ldots = \A_N$ by typing \verb|A=D2+2*diag(x)*D1+((2*N+1)/N)*eye(M);|. Note that \verb|D2*U| corresponds to $\partial^2u/\partial x_j^2$ in \eqref{e:uxt}; \verb|2*diag(x)*D1*U|, to $2x_j\partial u/\partial x_j$; and \verb|((2*N+1)/N)*eye(M)*U|, to $((2N+1)/N)u$, because we have divided the term $(2N+1)u$ in $N$ equal pieces, and assigned each piece to a matrix $\A_j$. With respect to \verb|B|, that corresponds to the term $-\exp(-\xx\cdot\xx)$, the simplest option is to create an $N$-dimensional mesh, which can be generated in MATLAB by typing \verb|XX = cell(1, N);|, followed by \verb|[XX{:}] = ndgrid(x);|. This syntax is especially useful, since it allows us to consider different values of $N$, without modifying the code. Then, \verb|XX{1}| corresponds to $x_1$, \verb|XX{2}|, to $x_2$, and so on, and it is straightforward to generate \verb|B|, after which \verb|XX|, which requires $N$ times as much storage as \verb|U| and is no longer used, can be removed from memory. The resulting code is
\begin{verbatim}
XX=cell(1,N);
[XX{:}]=ndgrid(x);
B=-exp(-XX{1}.^2);
for m=2:N
    B=B.*exp(-XX{m}.^2);
end
clear XX
\end{verbatim}
Even if in some specific cases it can be necessary to create explicitly the mesh, this is unnecessary in this example, thanks to the possibilities offered by MATLAB's syntax to multiply arrays of different dimensions, and to the command \verb|permute|. Indeed, it is enough to type
\begin{verbatim}
B=-exp(-x.^2);
for m=2:N
    B=B.*permute(exp(-x.^2),[m 2:m-1 1 m+1:N]);
end
\end{verbatim}
which produces an identical version of \verb|B|, but with fewer memory requirements, so we have used this approach.

Bearing in mind all the previous arguments, it is straightforward to implement a MATLAB code, which is offered in Listing~\ref{code:advdiftest}. Note that we impose that the numerical solution is real, to get rid of infinitesimally small imaginary parts (of the order of $\mathcal O(10^{-18})$) coming from round-off that might appear. The solution at $t = 1$, taking $N = 6$ dimensions, took $25.73$ seconds to execute, and the error was $9.6811\times10^{-14}$, which, in our opinion, is quite remarkable.

\begin{lstlisting}[label=code:advdiftest, language=MATLAB, basicstyle=\scriptsize, caption = {MATLAB program \texttt{advdiftest.m}, that solves numerically the advection--diffusion problem \eqref{e:uxt}}]
clear
tic
t=1; % time at which the numerical solution will be computed
N=6; % number of dimensions $N$
M=16; % number of nodes at each dimension
% create $\mathbf x$ and the differentiation matrices
[x,DD]=herdif(M,2,1.4);
D1=DD(:,:,1); % first-order differentiation matrix
D2=DD(:,:,2); % second-order differentiation matrix
AA=cell(1,N); % stores the matrices $\mathbf A_j$
B=-exp(-x.^2); % creates $\mathbf B$
for m=2:N
    B=B.*permute(exp(-x.^2),[m 2:m-1 1 m+1:N]);
end
U0=-2*B; % initial data
Uexact=-(1+exp(t))*B; % exact solution at time $t$
A=D2+2*diag(x)*D1+((2*N+1)/N)*eye(M); % creates $\mathbf A_j$
for m=1:N
    AA{m}=A;
end
U=real(evolND(AA,B,U0,t)); % computes the solution at time $t$
norm(U(:)-Uexact(:),inf) % error
toc % elapsed time
\end{lstlisting}

\subsection*{Author Contribution}

Both authors contributed equally to the conceptualization and writing. The computational implementation was carried out by Francisco de la Hoz. Critical editing was performed by Carlota M. Cuesta. Both authors have read and agreed to the published version of the manuscript.

\subsection*{Funding}

This work was partially supported by the research group grant IT1615-22 funded by the Basque Government, and by the projects PID2021-126813NB-I00 and PID2024-158099NB-I00 funded by MICIU/AEI/10.13039/501100011033 and by ERDF, EU.

\subsection*{Data Availability}

The MATLAB codes used in this paper are available in the GitHub repository \cite{githubSchurND2026}.

\subsection*{Declarations}

\subsubsection*{Competing interests}

The authors declare no competing interests.



\begin{thebibliography}{19}
	\ifx \bisbn   \undefined \def \bisbn  #1{ISBN #1}\fi
	\ifx \binits  \undefined \def \binits#1{#1}\fi
	\ifx \bauthor  \undefined \def \bauthor#1{#1}\fi
	\ifx \batitle  \undefined \def \batitle#1{#1}\fi
	\ifx \bjtitle  \undefined \def \bjtitle#1{#1}\fi
	\ifx \bvolume  \undefined \def \bvolume#1{\textbf{#1}}\fi
	\ifx \byear  \undefined \def \byear#1{#1}\fi
	\ifx \bissue  \undefined \def \bissue#1{#1}\fi
	\ifx \bfpage  \undefined \def \bfpage#1{#1}\fi
	\ifx \blpage  \undefined \def \blpage #1{#1}\fi
	\ifx \burl  \undefined \def \burl#1{\textsf{#1}}\fi
	\ifx \doiurl  \undefined \def \doiurl#1{\url{https://doi.org/#1}}\fi
	\ifx \betal  \undefined \def \betal{\textit{et al.}}\fi
	\ifx \binstitute  \undefined \def \binstitute#1{#1}\fi
	\ifx \binstitutionaled  \undefined \def \binstitutionaled#1{#1}\fi
	\ifx \bctitle  \undefined \def \bctitle#1{#1}\fi
	\ifx \beditor  \undefined \def \beditor#1{#1}\fi
	\ifx \bpublisher  \undefined \def \bpublisher#1{#1}\fi
	\ifx \bbtitle  \undefined \def \bbtitle#1{#1}\fi
	\ifx \bedition  \undefined \def \bedition#1{#1}\fi
	\ifx \bseriesno  \undefined \def \bseriesno#1{#1}\fi
	\ifx \blocation  \undefined \def \blocation#1{#1}\fi
	\ifx \bsertitle  \undefined \def \bsertitle#1{#1}\fi
	\ifx \bsnm \undefined \def \bsnm#1{#1}\fi
	\ifx \bsuffix \undefined \def \bsuffix#1{#1}\fi
	\ifx \bparticle \undefined \def \bparticle#1{#1}\fi
	\ifx \barticle \undefined \def \barticle#1{#1}\fi
	\bibcommenthead
	\ifx \bconfdate \undefined \def \bconfdate #1{#1}\fi
	\ifx \botherref \undefined \def \botherref #1{#1}\fi
	\ifx \url \undefined \def \url#1{\textsf{#1}}\fi
	\ifx \bchapter \undefined \def \bchapter#1{#1}\fi
	\ifx \bbook \undefined \def \bbook#1{#1}\fi
	\ifx \bcomment \undefined \def \bcomment#1{#1}\fi
	\ifx \oauthor \undefined \def \oauthor#1{#1}\fi
	\ifx \citeauthoryear \undefined \def \citeauthoryear#1{#1}\fi
	\ifx \endbibitem  \undefined \def \endbibitem {}\fi
	\ifx \bconflocation  \undefined \def \bconflocation#1{#1}\fi
	\ifx \arxivurl  \undefined \def \arxivurl#1{\textsf{#1}}\fi
	\csname PreBibitemsHook\endcsname
	
	%%% 1
	\bibitem[\protect\citeauthoryear{Sylvester}{1884}]{sylvester1884}
	\begin{barticle}
		\bauthor{\bsnm{Sylvester}, \binits{J.J.}}:
		\batitle{{Sur l'\'equation en matrices $px = xq$}}.
		\bjtitle{Comptes Rendus de l'Acad\'emie des Sciences}
		\bvolume{99}(\bissue{2}),
		\bfpage{67}--\blpage{71}
		(\byear{1884}).
		\bcomment{In French. Continued on pp. 115--116.}
	\end{barticle}
	\endbibitem
	
	%%% 2
	\bibitem[\protect\citeauthoryear{Bhatia and Rosenthal}{1997}]{bhatia1997}
	\begin{barticle}
		\bauthor{\bsnm{Bhatia}, \binits{R.}},
		\bauthor{\bsnm{Rosenthal}, \binits{P.}}:
		\batitle{{How and Why to Solve the Operator Equation $AX-XB = Y$}}.
		\bjtitle{Bulletin of the London Mathematical Society}
		\bvolume{29}(\bissue{1}),
		\bfpage{1}--\blpage{21}
		(\byear{1997})
	\end{barticle}
	\endbibitem
	
	%%% 3
	\bibitem[\protect\citeauthoryear{Li et~al.}{2010}]{benwen2010}
	\begin{barticle}
		\bauthor{\bsnm{Li}, \binits{B.-W.}},
		\bauthor{\bsnm{Tian}, \binits{S.}},
		\bauthor{\bsnm{Sun}, \binits{Y.-S.}},
		\bauthor{\bsnm{Hu}, \binits{Z.-M.}}:
		\batitle{{Schur-decomposition for 3D matrix equations and its application in
				solving radiative discrete ordinates equations discretized by Chebyshev
				collocation spectral method}}.
		\bjtitle{Journal of Computational Physics}
		\bvolume{229}(\bissue{4}),
		\bfpage{1198}--\blpage{1212}
		(\byear{2010})
	\end{barticle}
	\endbibitem
	
	%%% 4
	\bibitem[\protect\citeauthoryear{Bartels and Stewart}{1972}]{bartels1972}
	\begin{barticle}
		\bauthor{\bsnm{Bartels}, \binits{R.H.}},
		\bauthor{\bsnm{Stewart}, \binits{G.W.}}:
		\batitle{{Solution of the matrix equation $AX + XB = C$}}.
		\bjtitle{Communications of the ACM}
		\bvolume{15}(\bissue{9}),
		\bfpage{820}--\blpage{826}
		(\byear{1972})
	\end{barticle}
	\endbibitem
	
	%%% 5
	\bibitem[\protect\citeauthoryear{Higham}{2008}]{higham2008}
	\begin{bbook}
		\bauthor{\bsnm{Higham}, \binits{N.J.}}:
		\bbtitle{{F}unctions of {M}atrices. {T}heory and {C}omputation}.
		\bpublisher{Society for Industrial and Applied Mathematics},
		\blocation{Philadelphia}
		(\byear{2008})
	\end{bbook}
	\endbibitem
	
	%%% 6
	\bibitem[\protect\citeauthoryear{de~la Hoz and
		Vadillo}{2013}]{delahozvadillo2013a}
	\begin{barticle}
		\bauthor{\bsnm{Hoz}, \binits{F.}},
		\bauthor{\bsnm{Vadillo}, \binits{F.}}:
		\batitle{{A Sylvester-Based IMEX Method via Differentiation Matrices for
				Solving Nonlinear Parabolic Equations}}.
		\bjtitle{Communications in Computational Physics}
		\bvolume{14}(\bissue{4}),
		\bfpage{1001}--\blpage{1026}
		(\byear{2013})
	\end{barticle}
	\endbibitem
	
	%%% 7
	\bibitem[\protect\citeauthoryear{Stewart}{1996}]{stewart1996}
	\begin{botherref}
		\oauthor{\bsnm{Stewart}, \binits{G.W.}}:
		{Stochastic Automata, Tensors Operation, and Matrix Equations}.
		Technical Report, University of Maryland, UMIACS TR-96-11, CMSC TR-3598
		(1996)
	\end{botherref}
	\endbibitem
	
	%%% 8
	\bibitem[\protect\citeauthoryear{Touzene}{2006}]{touzene2006}
	\begin{bchapter}
		\bauthor{\bsnm{Touzene}, \binits{A.}}:
		\bctitle{{Approximated Tensor Sum Preconditioner for Stochastic Automata
				Networks}}.
		In: \bbtitle{IPDPS'06: Proceedings of the 20th International Conference on
			Parallel and Distributed Processing}
		(\byear{2006})
	\end{bchapter}
	\endbibitem
	
	%%% 9
	\bibitem[\protect\citeauthoryear{Zhang and Wang}{2021}]{xinfang2021}
	\begin{barticle}
		\bauthor{\bsnm{Zhang}, \binits{X.-F.}},
		\bauthor{\bsnm{Wang}, \binits{Q.-W.}}:
		\batitle{{Developing iterative algorithms to solve Sylvester tensor
				equations}}.
		\bjtitle{Applied Mathematics and Computation}
		\bvolume{409},
		\bfpage{126403}
		(\byear{2021})
	\end{barticle}
	\endbibitem
	
	%%% 10
	\bibitem[\protect\citeauthoryear{Chen and Li}{2024}]{YuhanChen2023}
	\begin{barticle}
		\bauthor{\bsnm{Chen}, \binits{Y.}},
		\bauthor{\bsnm{Li}, \binits{C.}}:
		\batitle{{A Tensor Multigrid Method for Solving Sylvester Tensor Equations}}.
		\bjtitle{IEEE Transactions on Automation Science and Engineering}
		\bvolume{21}(\bissue{3}),
		\bfpage{4397}--\blpage{4405}
		(\byear{2024})
	\end{barticle}
	\endbibitem
	
	%%% 11
	\bibitem[\protect\citeauthoryear{Chen and Kressner}{2020}]{chen2020}
	\begin{barticle}
		\bauthor{\bsnm{Chen}, \binits{M.}},
		\bauthor{\bsnm{Kressner}, \binits{D.}}:
		\batitle{{Recursive blocked algorithms for linear systems with Kronecker
				product structure}}.
		\bjtitle{Numerical Algorithms}
		\bvolume{84},
		\bfpage{1199}--\blpage{1216}
		(\byear{2020}).
		\bcomment{The codes are available at
			\url{https://www.epfl.ch/labs/anchp/wp-content/uploads/2019/05/tensor_recursive.tar.gz}}
	\end{barticle}
	\endbibitem
	
	%%% 12
	\bibitem[\protect\citeauthoryear{Massei and Robol}{2024}]{massei2024}
	\begin{barticle}
		\bauthor{\bsnm{Massei}, \binits{S.}},
		\bauthor{\bsnm{Robol}, \binits{L.}}:
		\batitle{{A nested divide-and-conquer method for tensor Sylvester equations
				with positive definite hierarchically semiseparable coefficients}}.
		\bjtitle{IMA Journal of Numerical Analysis}
		\bvolume{44}(\bissue{6}),
		\bfpage{3482}--\blpage{3519}
		(\byear{2024})
	\end{barticle}
	\endbibitem
	
	%%% 13
	\bibitem[\protect\citeauthoryear{Ballani and Grasedyck}{2013}]{Ballani2013}
	\begin{barticle}
		\bauthor{\bsnm{Ballani}, \binits{J.}},
		\bauthor{\bsnm{Grasedyck}, \binits{L.}}:
		\batitle{{A projection method to solve linear systems in tensor format}}.
		\bjtitle{Numerical linear algebra with applications}
		\bvolume{20}(\bissue{1}),
		\bfpage{27}--\blpage{43}
		(\byear{2013})
	\end{barticle}
	\endbibitem
	
	%%% 14
	\bibitem[\protect\citeauthoryear{de~la Hoz and
		Vadillo}{2013}]{delahozvadillo2013b}
	\begin{barticle}
		\bauthor{\bsnm{Hoz}, \binits{F.}},
		\bauthor{\bsnm{Vadillo}, \binits{F.}}:
		\batitle{{The solution of two-dimensional advection-diffusion equations via
				operational matrices}}.
		\bjtitle{Applied Numerical Mathematics}
		\bvolume{72},
		\bfpage{172}--\blpage{187}
		(\byear{2013})
	\end{barticle}
	\endbibitem
	
	%%% 15
	\bibitem[\protect\citeauthoryear{Grasedyck}{2004}]{Grasedyck2004}
	\begin{barticle}
		\bauthor{\bsnm{Grasedyck}, \binits{L.}}:
		\batitle{{Existence and Computation of Low Kronecker-Rank Approximations for
				Large Linear Systems of Tensor Product Structure}}.
		\bjtitle{Computing}
		\bvolume{72},
		\bfpage{247}--\blpage{265}
		(\byear{2004})
	\end{barticle}
	\endbibitem
	
	%%% 16
	\bibitem[\protect\citeauthoryear{Cuesta and de~la
		Hoz}{2026}]{githubSchurND2026}
	\begin{botherref}
		\oauthor{\bsnm{Cuesta}, \binits{C.M.}},
		\oauthor{\bsnm{Hoz}, \binits{F.}}:
		{SchurND}.
		GitHub repository.
		\url{https://github.com/fdlhm/SchurND}
		(2026)
	\end{botherref}
	\endbibitem
	
	%%% 17
	\bibitem[\protect\citeauthoryear{Knuth}{2011}]{Knuth2011}
	\begin{bbook}
		\bauthor{\bsnm{Knuth}, \binits{D.E.}}:
		\bbtitle{The Art of Computer Programming, Volume~4A: Combinatorial Algorithms,
			Part~1}.
		\bpublisher{Addison-Wesley},
		\blocation{Upper Saddle River, NJ}
		(\byear{2011})
	\end{bbook}
	\endbibitem
	
	%%% 18
	\bibitem[\protect\citeauthoryear{Parlett}{1976}]{parlett1976}
	\begin{barticle}
		\bauthor{\bsnm{Parlett}, \binits{B.N.}}:
		\batitle{{A recurrence among the elements of functions of triangular
				matrices}}.
		\bjtitle{Linear Algebra and its Applications}
		\bvolume{14}(\bissue{2}),
		\bfpage{117}--\blpage{121}
		(\byear{1976})
	\end{barticle}
	\endbibitem
	
	%%% 19
	\bibitem[\protect\citeauthoryear{Weideman and Reddy}{2000}]{WeidemanReddy2000}
	\begin{barticle}
		\bauthor{\bsnm{Weideman}, \binits{J.A.C.}},
		\bauthor{\bsnm{Reddy}, \binits{S.C.}}:
		\batitle{{A MATLAB differentiation matrix suite}}.
		\bjtitle{ACM Transactions on Mathematical Software}
		\bvolume{26}(\bissue{4}),
		\bfpage{465}--\blpage{519}
		(\byear{2000}).
		\bcomment{The codes are available at
			\url{https://appliedmaths.sun.ac.za/~weideman/research/differ.html}}
	\end{barticle}
	\endbibitem
	
\end{thebibliography}
\end{document}